\documentclass[10pt]{amsart}
\usepackage[margin=1.5in]{geometry}
\usepackage{amsmath,amssymb,amsthm}
\usepackage[dvipsnames]{xcolor}
\usepackage[colorlinks, allcolors=purple!90!black]{hyperref} 
\hypersetup{
    colorlinks = true,
    linkcolor = Blue,
    citecolor = purple!90!black,
    urlcolor = OliveGreen,
}

\usepackage[parfill]{parskip}

\makeatletter
\renewcommand{\@cite}[2]{{\m@th\upshape\bfseries[#1\if\relax\detokenize{#2}\relax\else, #2\fi]}}
\makeatother
\newtheorem{theorem}{Theorem}
\newtheorem{proposition}{Proposition}
\newtheorem{lemma}{Lemma}
\numberwithin{equation}{section}
\begin{document}
\title[Gaps between zeros of Dirichlet $L$-functions]{Small and large gaps between zeros of Dirichlet $L$-functions}
\author[B. Durkan]{Benjamin Durkan}
\address{Department of Mathematics, The University of Manchester, Oxford Road, Manchester, M13 9PL, United Kingdom}
\email{benjamin.durkan@postgrad.manchester.ac.uk}
\author[A. Pearce-Crump]{Andrew Pearce-Crump}
\address{School of Mathematics, University of Bristol, Fry Building, Woodland Road, Bristol, BS8 1UG, United Kingdom}
\email{andrew.pearce-crump@bristol.ac.uk}
\date{}
\begin{abstract}
For every fixed primitive Dirichlet $L$-function and every $\varepsilon>0$, we prove unconditionally that infinitely many consecutive zeros have spacing less than $1/2+\varepsilon$ times the mean. Assuming the Riemann Hypothesis, we also prove a positive proportion of such spacings for the Riemann zeta-function, and infinitely many spacings greater than $3.63846723$ times the mean for every fixed primitive Dirichlet $L$-function.
\end{abstract}
\maketitle
\section{Introduction and main results}\label{sec:results}

The distribution of zeros of $L$-functions is a central problem in analytic number theory. In this paper we study gaps between consecutive zeros of a fixed primitive Dirichlet $L$-function. Let $\chi$ have conductor $Q$, with $Q=1$ corresponding to the Riemann zeta-function, and write $0<\gamma_1\le\gamma_2\le\cdots$ for the positive ordinates of the non-trivial zeros of $L(s,\chi)$, counted with multiplicity. Set
$$
\mu=\liminf_{n\to\infty}\frac{(\gamma_{n+1}-\gamma_n)\log\gamma_n}{2\pi},
\qquad
\lambda=\limsup_{n\to\infty}\frac{(\gamma_{n+1}-\gamma_n)\log\gamma_n}{2\pi}.
$$
By the zero-counting formula \cite{Davenport,Titchmarsh}, the mean spacing at height $T$ is asymptotic to $2\pi/\log T$, since $Q$ is fixed. Our first result reaches half this spacing.

\begin{theorem}\label{thm:small}
Unconditionally, $\mu\le1/2$ for every fixed primitive Dirichlet character $\chi$.
\end{theorem}

By the functional equation, zeros off the critical line occur in distinct pairs with the same ordinate, even for complex characters. If there are infinitely many such pairs at positive heights, then, since we count multiplicities, there are infinitely many zero spacings and hence $\mu=0$. We may therefore assume in proving Theorem~\ref{thm:small} that all zeros above some fixed positive height lie on the critical line.

Under the Riemann Hypothesis, the small-gap result of Theorem \ref{thm:small} for zeta can be strengthened to a positive-proportion result.

\begin{theorem}\label{thm:positive-proportion}
Assume the Riemann Hypothesis for the Riemann zeta-function. For every $\varepsilon>0$, there is a constant $c_\varepsilon>0$ such that, for all sufficiently large $T$,
$$
\#\left\{n:T\le\gamma_n\le\gamma_{n+1}\le2T,\quad
\frac{(\gamma_{n+1}-\gamma_n)\log\gamma_n}{2\pi}<\frac12+\varepsilon\right\}
\ge c_\varepsilon T\log T,
$$
where the $\gamma_n$ are the positive non-trivial zero ordinates of zeta, counted with multiplicity.
\end{theorem}

The count includes spacings of length zero, and $c_\varepsilon$ may tend to zero as $\varepsilon\to0$. Thus the theorem does not give a positive proportion of gaps at most half the mean spacing.

Our final theorem shows, conditional on the Riemann Hypothesis, that there exist infinitely many gaps which are much larger than the average spacing.

\begin{theorem}\label{thm:large}
Assume the Riemann Hypothesis for $L(s,\chi)$. Then $\lambda>3.63846723$.
\end{theorem}

For $\zeta(s)$, Selberg \cite{Selberg} stated, without assuming the Riemann Hypothesis, that $\mu<1<\lambda$. Under the Riemann Hypothesis, Montgomery and Odlyzko \cite{MO} obtained $\mu<0.5179$ and $\lambda>1.9799$. Conrey, Ghosh and Gonek \cite{CGG84} simplified their method and obtained $\mu<0.5172$ and $\lambda>2.337$. Further choices of arithmetic coefficients led to the bounds in \cite{BMN,FW12,Preobrazhenskii}. Inoue \cite{Inoue} combined the method with pair correlation and obtained $\mu<0.50895$.

For large gaps, comparisons of discrete and continuous mean values were developed in \cite{Mueller,CGG86,Ng,Bui11,FW13,Bui14}, with the results of \cite{CGG86,Ng,Bui11,FW13} assuming the Generalised Riemann Hypothesis. Hall \cite{Hall}, Bredberg \cite{Bredberg}, and Bui and Milinovich \cite{BM} used Wirtinger's inequality to study consecutive zeros on the critical line, unconditionally. The last of these obtained a bound of $3.18$, which applies to $\lambda$ under the Riemann Hypothesis. For explicit unconditional bounds, see \cite{STTB}. Positive proportions and gaps between distinct zeros are treated in \cite{CGGGH,Sou96,Wu,BGMM}, and gaps between zeros that are $r$ apart in \cite{CTB,Inoue24,IKT}.

The classical Montgomery--Odlyzko criterion uses a weighted mean of the number of zeros in a short interval \cite{MO}. Goldston, Trudgian and Turnage-Butterbaugh \cite{GTT} proved that, for polynomials of length $T^{1-\delta}$ with fixed $\delta>0$, no choice of coefficients allows this criterion to give a small-gap bound below $0.5042$. For polynomials of length at most $T$, Inoue, Kobayashi and Toma \cite{IKT} strengthened the small-gap restriction to $0.508$ and obtained a large-gap restriction at $3.022$. Our small-gap proof uses all fixed moments of the logarithmic difference. For large gaps we use sums over reduced fractions $m/n$, so multiplication by a prime either adds a factor to the numerator or cancels one from the denominator, and symmetric coefficients make the two contributions equal. 

Bondarenko and Heap \cite{BondarenkoHeapHalf} recently proved independently that $\mu\le1/2$ for $\zeta(s)$ under the Riemann Hypothesis by applying resonance to a quotient of shifted zeta-functions. If infinitely many zeros lie off the critical line at positive heights, their paired ordinates give $\mu=0$. Otherwise, all sufficiently high positive zeros lie on the critical line, and their proof applies. Here we compute every fixed moment of a shifted logarithmic difference, using the squared modulus of a M\"obius Dirichlet polynomial as the weight. If $\mu>1/2$, these must be the moments of a positive measure on two horizontal lines. By Lemma~\ref{lem:two-lines}, however, the moments determine a density which is negative on part of the lower line, giving a contradiction. This calculation also applies to every fixed primitive Dirichlet character.

Montgomery \cite{Montgomery} proved under the Riemann Hypothesis that the pair correlation of zeta zeros agrees with the GUE prediction for test functions with restricted Fourier support. The full GUE hypothesis predicts the same local statistics as the eigenvalues of large random Hermitian matrices from the Gaussian unitary ensemble \cite{LagariasRodgers}. In particular, it predicts a positive proportion of consecutive gaps below every fixed positive multiple of the mean spacing and above every fixed multiple, giving $\mu=0$ and $\lambda=\infty$.

Under the Alternative Hypothesis, consecutive normalised spacings approach positive half-integers \cite{LagariasRodgers}. In the formulation used by Baluyot, who assumes the Riemann Hypothesis and counts multiplicities, zero is also allowed \cite{BaluyotAH}. Lagarias and Rodgers construct a half-integer model satisfying the restricted correlation formulas considered in their paper. Theorems~\ref{thm:small}--\ref{thm:large} remain compatible with this hypothesis. The small-gap bound reaches the first positive half-integer, while Theorem~\ref{thm:positive-proportion} allows the counted gaps to approach $1/2$ from above. The large-gap bound allows spacings near larger half-integers.

\section{Overview of the proof}\label{sec:sketch}

Suppose, contrary to Theorem~\ref{thm:small}, that $\mu>1/2$, and choose
$$
\frac12<\phi<\min(\mu,1),\qquad h=\frac{2\pi\phi}{\log T}.
$$
All sufficiently high positive zeros then lie on the critical line, and every interval of length $h$ in $[T-1,2T+1]$ contains at most one zero. We will express this count as the imaginary part of a logarithmic difference and calculate all its fixed moments. These moments would have to come from a positive measure on the lines $\operatorname{Im}w=0,1$. We show that they determine a density which is negative on part of one line, contradicting the gap assumption.

Let $\xi_\chi$ denote the entire completed $L$-function, with the constant factor chosen so that it is real on the critical line, and define
$$
V_\phi(t)=\frac1\pi\left(\log\xi_\chi\left(\frac12+i(t+h/2)\right)
-\log\xi_\chi\left(\frac12+i(t-h/2)\right)\right).
$$
We choose this logarithmic difference because its imaginary part counts the zeros in $(t-h/2,t+h/2)$, with multiplicity. Here the logarithms are boundary values from $\operatorname{Re}s>1/2$, and neither endpoint is a zero ordinate. Under our gap assumption, the interval contains at most one zero, so $V_\phi(t)$ lies on one of the two lines $\operatorname{Im}w=0,1$. The advantage of using the full complex difference is that we can calculate its moments by contour integration.

The weight is the squared modulus of
$$
R_T(t)=\sum_{\substack{n\le N\\(n,Q)=1}}
\frac{\chi(n)\mu(n)}{\sqrt n}
f\left(\frac{\log n}{\log N}\right)n^{-it},
\qquad N=\left\lfloor\frac{T}{(\log T)^{20}}\right\rfloor,
$$
where $\mu(n)$ is the M\"obius function and $f$ is a nonzero continuous function on $[0,1]$, to be chosen later. On shifting to $\operatorname{Im}t=-1$, we may expand the logarithmic difference in its absolutely convergent Euler series. Integration with the weight in \eqref{eq:gaussian-weight} leaves the terms with $n'=kn$ as the main terms. For the other terms, Proposition~\ref{prop:finite-moments} uses a lower bound for $|\log(kn/n')|$ to prove an error estimate uniform in the coefficients.

The M\"obius coefficients restrict the initial and final integers to be squarefree, so only multiplication by distinct primes contributes. For $p\nmid n$, the identity $\mu(np)=-\mu(n)$ reverses the sign of a prime term. With $u=\log p/\log N$, Mertens' estimate changes prime summation into integration with respect to $du/u$. For repeated multiplication we therefore apply the same integral repeatedly. Write
$$
(A_\phi f)(x)=i\phi f(x)+\frac{2i}\pi\int_0^x
\frac{\sin(\pi\phi u)}u f(x-u)du,\qquad 0\le x\le1.
$$
By Proposition~\ref{prop:integer-limit}, the $j$th weighted moment tends to
$$
\frac{\int_0^1(A_\phi^jf)(x)\overline{f(x)}dx}{\int_0^1|f(x)|^2dx}.
$$
The first term in $A_\phi$ comes from the completion factor, and the integral from multiplication by a prime. We write $A_\phi^j$ for $j$ successive applications of this formula.

Bounds for the absolute moments allow us to pass to a positive limiting measure with these moments. To recover its density, we use $w\mapsto e^{\pi w}$, which sends the lower line to the positive real axis and the upper line to the negative real axis. Lemma~\ref{lem:two-lines} uses Cauchy's formula after this substitution and proves that the measure is uniquely determined.

It remains to choose $f$ so that the density is negative. The exponential of $A_\phi$, defined by its convergent power series, is particularly simple on $e^{-i\pi\phi x}$. By \eqref{eq:half-exponential}, it sends this function to $e^{i\pi\phi}e^{i\pi\phi x}$. Taking their difference gives
$$
f(x)=e^{-i\pi\phi x}-e^{i\pi\phi}e^{i\pi\phi x}
=-2ie^{i\pi\phi/2}\sin\left(\pi\phi(x+1/2)\right).
$$
This choice cancels the inverse occurring in the density formula \eqref{eq:two-line-density}. At zero on the lower line we obtain
\begin{equation}\label{eq:outline-density}
p_0(0)=\frac{\sin(\pi\phi)\sin(2\pi\phi)}
{\pi\phi\displaystyle\int_0^1|f(x)|^2dx}<0.
\end{equation}
For $1/2<\phi<1$, the density is negative on an interval about zero, which is impossible for a positive measure. The half-spacing bound comes from the change of sign of $\sin(2\pi\phi)$.

For Theorem~\ref{thm:positive-proportion}, we assume the Riemann Hypothesis for zeta and choose $1/2<\phi<\min(1,1/2+\varepsilon)$. We need to show that the contradiction comes from many short intervals, rather than a set of centres of very small measure. We shorten the polynomial to $T^\theta$, where $1/\phi-1<\theta<1$, and remove prime factors below $T^\eta$, for a small fixed $\eta>0$. Each remaining integer has a bounded number of prime factors. The explicit formula and local zero-count estimates then bound every fixed absolute moment of both the polynomial and $V_\phi$.

For fixed $\eta$, we first let $T\to\infty$. Letting $\eta\downarrow0$ afterwards recovers the moment formulas on $[0,\theta]$. The numerator in \eqref{eq:outline-density} becomes $\sin(\pi\phi\theta)\sin(\pi\phi(1+\theta))$, which is negative by the choice of $\theta$. We use this density to construct a polynomial $P$ whose real part is nonnegative on both lines but whose calculated weighted mean has negative real part. We then fix $P$, a smooth coefficient function and a sufficiently small $\eta$. By H\"older's inequality and the absolute-moment bounds, intervals containing at least two zeros have centres of measure $\gg T$. A fixed consecutive pair belongs to such intervals for centres of measure at most $h$, giving $\gg T\log T$ small spacings.

For Theorem~\ref{thm:large}, a weighted mean of the number of zeros below one is sufficient, since it forces an empty interval. Assume the Riemann Hypothesis and use reduced fractions $m/n$ with squarefree coprime $m,n\le Y=\sqrt T/(\log T)^{10}$. The coefficients are symmetric in $m,n$. They depend on the logarithmic sizes of these integers and on four sums describing how their sizes are divided among their prime factors, as in \eqref{eq:large-marked-coefficients}. Two explicit quadratic forms calculate the mean count. The numerical inequality \eqref{eq:marked-large-margin} makes this mean less than one at $\phi=3.63846723$. By continuity the inequality holds at a fixed larger $\phi$, proving the strict bound in Theorem~\ref{thm:large}.

We begin with the boundary values of the logarithms and the required harmonic sums in Section~\ref{sec:preliminaries}. Sections~\ref{sec:moments} and~\ref{sec:integer-moments} evaluate the moments by contour integration and prime summation. Section~\ref{sec:two-line-moments} recovers the density and completes the proof of Theorem~\ref{thm:small}. The additional estimates for Theorem~\ref{thm:positive-proportion} are proved in Section~\ref{sec:positive-proportion}. Section~\ref{sec:large} contains the large-gap calculation and the proof of Theorem~\ref{thm:large}.

\section{Preliminary lemmas}\label{sec:preliminaries}

Fix a primitive character $\chi$ of conductor $Q$, including $Q=1$, and write
$$
\xi_\chi(s)=\varepsilon_\chi^{-1/2}\left(\frac Q\pi\right)^{(s+a_\chi)/2}
\Gamma\left(\frac{s+a_\chi}{2}\right)L(s,\chi),
$$
where $\chi(-1)=(-1)^{a_\chi}$, $a_\chi\in\{0,1\}$, and $\varepsilon_\chi$ is the root number \cite{Davenport}. Choose either square root of $\varepsilon_\chi$. When $Q=1$, include the factor $s(s-1)/2$ to obtain the Riemann $\xi$-function. Since $\varepsilon_\chi^{-1/2}\varepsilon_\chi=\overline{\varepsilon_\chi^{-1/2}}$, the functional equation and $L(s,\overline\chi)=\overline{L(\overline s,\chi)}$ imply $\xi_\chi(s)=\overline{\xi_\chi(1-\overline s)}$. Thus $\xi_\chi$ is real on the critical line, also when $\chi$ is complex. The notation $\boldsymbol1$ denotes one or zero according as its subscripted condition holds or fails.

\begin{lemma}\label{lem:logarithmic-boundary}
Suppose all zeros above a fixed height $Y_0\ge2$ lie on the critical line. Continue the Euler-series and completion logarithms into the simply connected zero-free region $\operatorname{Re}s>1/2$, $\operatorname{Im}s>Y_0$. For fixed $\phi>0$, set $h=2\pi\phi/\log T$ and
$$
V_\phi(t)=\frac1\pi\left(\log\xi_\chi\left(\frac12+i(t+h/2)\right)
-\log\xi_\chi\left(\frac12+i(t-h/2)\right)\right).
$$
Boundary values from $\operatorname{Re}s>1/2$ exist almost everywhere. For $t\asymp T$, $\operatorname{Im}V_\phi(t)$ counts zeros with multiplicity in $(t-h/2,t+h/2)$ unless an endpoint is a zero. For $u\asymp T$ and every fixed integer $j\ge1$,
$$
\sup_{0<\epsilon\le1}\int_u^{u+1}|\log\xi_\chi(1/2+\epsilon+it)|^jdt
\ll_{\chi,Y_0,j}T^{j+1}.
$$
\end{lemma}
\begin{proof}
We first bound the logarithm near its boundary zeros. Set $X(z)=\xi_\chi(1/2+iz)$. The functional equation implies $X(\overline z)=\overline{X(z)}$, and $X$ is an entire function of order one. By Hadamard's factorisation theorem \cite[Theorem~5.6]{IK}\cite[Theorem~2.12 and (2.12.5)]{Titchmarsh}, we have
$$
X(z)=e^{A_0+B_0z}z^m\prod_{\rho\ne0}(1-z/\rho)e^{z/\rho}.
$$
Here $m$ is the order at zero, and the remaining zeros $\rho$ are counted with multiplicity, with nonreal zeros grouped in conjugate pairs. We have $|\operatorname{Im}\rho|\le1/2$. There are $O_\chi(\log(2+|v|))$ zeros in a unit interval at height $v$. This follows from the zero-counting formulas in \cite{BennettZeros} when $Q>1$, and \cite[Theorems~9.2--9.4]{Titchmarsh} when $Q=1$, by taking the difference of the counts at $|v|+2$ and $|v|-2$. Logarithmic differentiation at a real nonzero point away from the zeros shows that $B_0$ is real. On the real axis, the product for a nonreal conjugate pair is $|1-t/\rho|^2e^{2t\operatorname{Re}(1/\rho)}>0$, so its argument is constant.

For $t\in[u,u+1]$, $u\asymp T$ and $0<\epsilon\le1$, the product bounds the chosen logarithmic branch by
$$
\left|\log\xi_\chi(1/2+\epsilon+it)\right|
\le C_{\chi,Y_0}T\log^2T+
C_\chi\sum_{\substack{\gamma\in[u-1,u+2]\\X(\gamma)=0}}|\log|t-\gamma||.
$$
For zeros with $|\rho|>6T$, the logarithm of the factor $(1-z/\rho)e^{z/\rho}$ is $O(T^2/|\rho|^2)$, and by the zero count we have
$$
T^2\sum_{|\rho|>6T}|\rho|^{-2}=O_\chi(T\log T).
$$
For the remaining factors away from local real zeros, the sum of $|\rho|^{-1}$ over $0<|\rho|\le6T$ is $O_\chi(\log^2T)$, giving $O_\chi(T\log^2T)$. Every nonreal zero of $X$ has real part at most $Y_0$ and is therefore at distance $\gg T$ from $[u,u+1]$. Near a real zero $\gamma$, set $d=t-\gamma$. We have $|\operatorname{Log}(\epsilon+id)|\le|\log|d||+\pi/2+(\log2)/2$. The value of the logarithmic branch at $\epsilon=1$ is $O_\chi(T\log T)$ by the Euler series and Stirling's formula.

There are $O_\chi(\log T)$ local zeros and $\int_{-2}^2|\log|v||^j\,dv<\infty$. Applying $\left(\sum_1^M z_\nu\right)^j\le M^{j-1}\sum_1^M z_\nu^j$ gives an integrable majorant independent of $\epsilon$, and hence
$$
\sup_{0<\epsilon\le1}\int_u^{u+1}|\log\xi_\chi(1/2+\epsilon+it)|^jdt
\ll_{\chi,Y_0,j}T^{j+1}.
$$
By dominated convergence, we may pass to the boundary in integrals of fixed powers at $t\pm h/2$. The real exponential factors and nonreal pairs have constant argument, and each intervening real zero increases the boundary argument from below by $\pi$, counted with multiplicity. Thus $\operatorname{Im}V_\phi$ counts the zeros unless an endpoint is a zero.

\end{proof}

\begin{lemma}\label{lem:arithmetic-sums}
Let $f,F$ be continuous complex functions on $[0,1]$. Then for fixed $Q$ we have
$$
\begin{aligned}
\lim_{N\to\infty}\frac1{\log N}\sum_{\substack{n\le N^{1-v}\\(n,Q)=1}}
\frac{\mu(n)^2}{n}
f\left(\frac{\log n}{\log N}\right)
\overline{f\left(\frac{\log n}{\log N}+v\right)}
=
\frac6{\pi^2}\prod_{p\mid Q}\frac p{p+1}
\int_0^{1-v}f(x)\overline{f(x+v)}\,dx
\end{aligned}
$$
uniformly for $0\le v\le1$, with both sides extended by zero for $v>1$. Also,
$$
\lim_{N\to\infty}\frac1{\log N}\sum_{\substack{p\le N\\p\nmid Q}}
\frac{\log p}{p}F\left(\frac{\log p}{\log N}\right)
=\int_0^1F(u)\,du.
$$
\end{lemma}
\begin{proof}
We use the classical estimates
$$
\sum_{\substack{n\le X\\(n,Q)=1}}\frac{\mu(n)^2}{n}
=\frac6{\pi^2}\prod_{p\mid Q}\frac p{p+1}\log X+O_Q(1),
\qquad
\sum_{\substack{p\le X\\p\nmid Q}}\frac{\log p}p=\log X+O_Q(1),
$$
which can be found in \cite{Davenport,Tenenbaum}, for example. The second is Mertens' estimate with the primes dividing $Q$ removed.

Taking $X=N^x$ and dividing the first estimate by $\log N$ gives an error $O_Q(1/\log N)$ uniformly for $0\le x\le1$. Approximate $f(x)\overline{f(x+v)}$ by step functions. Uniform continuity makes the approximation uniform in $v$, proving the first limit. At $v=1$, the only possible summand is $O_f(1/\log N)$, while the integral is zero. The second limit follows by partial summation of Mertens' estimate.
\end{proof}

\section{Beginning the proof}\label{sec:moments}

Assume that all zeros above a fixed height $Y_0\ge2$ lie on the critical line. We now evaluate the weighted moments of $V_\phi$ for fixed $\phi>0$. Here $h$ and $V_\phi$ are as in Lemma~\ref{lem:logarithmic-boundary}, so that $\operatorname{Im}V_\phi$ counts the zeros in $(t-h/2,t+h/2)$ almost everywhere.

We use integer indices for small gaps and reduced fractions for large gaps, treating both mean values together. For a reduced positive fraction $q=m/n$ with $(mn,Q)=1$, set $\chi(q)=\chi(m)\overline{\chi(n)}$, and let
$$
N=\left\lfloor\frac{T}{(\log T)^{20}}\right\rfloor,
\qquad Y=\frac{\sqrt T}{(\log T)^{10}},\qquad
R_T(t)=\sum_q\chi(q)c_qq^{-it}.
$$
The indices are either integers $q\le N$ with $(q,Q)=1$, or reduced fractions $m/n$ with $m,n\le Y$ and $(mn,Q)=1$. The coefficients $c=(c_q)$ may depend on $T$ and are not all zero. We denote by $\mathbb E_T$ the weighted average on $[T,2T]$ with weight
\begin{equation}\label{eq:gaussian-weight}
|R_T(t)|^2\exp\left(-\frac{(t-3T/2)^2(\log T)^2}{2T^2}\right)dt,
\end{equation}
normalised to have integral one.

Write $\Lambda(n)=\log p$ when $n=p^j$ for a prime $p$ and integer $j\ge1$, and $\Lambda(n)=0$ otherwise. On $\operatorname{Im}t=-1$ the logarithmic $L$-difference has the absolutely convergent series
$$
-2i\sum_{k\ge2}\frac{\chi(k)\Lambda(k)}{k^{1/2+it}\log k}
\sin\left(\frac{\pi\phi\log k}{\log T}\right).
$$
After multiplication by $|R_T(t)|^2$, the terms with no oscillation satisfy $r=kq$. Their character factors cancel, since $\chi(q)\overline{\chi(kq)}\chi(k)=1$. We record their coefficients in the matrix
$$
K_{r,q}=-\frac{2i\Lambda(k)}{\sqrt k\log k}
\sin\left(\frac{\pi\phi\log k}{\log T}\right)
\quad(r=kq,\ k\ge2),
$$
and set all other entries to zero. The powers $K^j$ sum over $j$ successive prime-power multiplications. All intermediate indices are included, even when $c_q=0$.

We write $\|c\|^2=\sum_q|c_q|^2$ and $\|f\|_2^2=\int_0^1|f(x)|^2\,dx$. The notation $i\phi+K/\pi$ means that $i\phi$ is added to every diagonal entry of $K/\pi$.

\begin{proposition}\label{prop:finite-moments}
For every fixed integer $j\ge0$,
$$
\mathbb E_T[V_\phi^j]
=\frac{\sum_{q,q'}\overline{c_{q'}}
\left(\left(i\phi+K/\pi\right)^j\right)_{q',q}c_q}
{\sum_q|c_q|^2}+o(1),
$$
uniformly over coefficient vectors which are not identically zero, on either index set. The error may depend on $\chi,Y_0,\phi,j$.

\end{proposition}
\begin{proof}
We show that only the terms $q'=kq$ contribute to the mean. Set $L=\log T$, $H=T/L$, $B=T/L^{20}$ and $g(t)=e^{-(t-3T/2)^2/(2H^2)}$. Either index set has at most $B$ elements and $q'/q\le B$. For fractions $q=m/n$ and $q'=m'/n'$, we have $mn',nm'\le Y^2=B$. For integers, $q,q'\le N\le B$. Hence the reduced denominator of $kq/q'$ is at most $B$ for every integer $k\ge1$. When $kq/q'\ne1$, its numerator and denominator are distinct positive integers, and therefore
$$
|\log(kq/q')|\ge\frac1{2B}.
$$

We first bound the sums of absolute values of the terms with $q'=kq$. Write $|K|$ for the matrix of absolute values and $\|M\|=\sup_{\|u\|=1}\|Mu\|$. When $r=p^e q$, $e\ge1$, the corresponding entry is at most
$$
\frac{2\pi\phi\log p}{Lp^{e/2}}.
$$
Write $r=m_r/n_r$ in lowest terms. Multiplication by $p^e$ either adds a numerator factor, cancels a denominator factor, or does both, changing $m/(p^bn)$ to $p^am/n$, where $a,b\ge1$, $e=a+b$ and $p\nmid mn$.

Let $x,y$ be nonnegative unit vectors. For numerator terms, by
$2x_rp^{-e/2}y_q\le x_r^2+p^{-e}y_q^2$, we have
$$
\begin{aligned}
\sum_{\mathrm{numerator\ terms}}|K_{r,q}|x_ry_q
&\le\frac{\pi\phi}L
\left(
\sum_rx_r^2\sum_{p^e\mid m_r}\log p+
\sum_qy_q^2\sum_{p^e\le T}\frac{\log p}{p^e}
\right)\\
&=O_\phi(1).
\end{aligned}
$$
For fixed $(p,e)$, specifying either $q$ or $r$ determines the other. The two sums are bounded by $\sum_{p^e\mid m_r}\log p=\log m_r\le L$ and $\sum_{p^e\le T}\log p/p^e=L+O(1)$. For denominator terms, interchange rows and columns and use $\sum_{p^e\mid n_q}\log p=\log n_q\le L$.

In the third case, for a fixed row, $a=v_p(m_r)\ge1$ is fixed, where $v_p(m_r)$ is the exponent of $p$ in $m_r$. By summing over all
$b\ge1$, we have
$$
\sum_{b\ge1}\frac{2\pi\phi\log p}{Lp^{(a+b)/2}}
\le\frac{C_\phi\log p}{Lp}.
$$
For a fixed column, hold $b=v_p(n_q)\ge1$ fixed and sum over $a\ge1$ for the same bound. Further size restrictions only decrease these sums. If $\log v\le L$, then
$$
\sum_{p\mid v}\frac{\log p}p
\le\sum_{p\le L}\frac{\log p}p+
\frac1L\sum_{p\mid v}\log p
\le\log L+O(1).
$$
Thus the row and column sums in the third case are $O_\phi(\log L/L)$. By Cauchy's inequality, its norm satisfies the same bound. Combining the cases, we have $\||K|\|=O_\phi(1)$.

The part of $|K|$ from proper prime powers, with $e\ge2$, has norm $O_\phi(L^{-1/2})$. To prove this, note that
$$
\sum_{\substack{p,\ e\ge2}}\frac{\log p}{p^e}
=\sum_p\frac{\log p}{p(p-1)}
$$
converges. For numerator terms use
$$
2x_rp^{-e/2}y_q
\le L^{-1/2}x_r^2+L^{1/2}p^{-e}y_q^2.
$$
The first resulting sum is $O(L^{1/2})$ by $\sum_{p^e\mid m_r}\log p=\log m_r\le L$, and the second is $O(L^{1/2})$ by the convergent prime-power sum. Multiplication by $O_\phi(L^{-1})$ gives $O_\phi(L^{-1/2})$. Interchanging rows and columns proves the denominator estimate. The third case already has norm $O_\phi(\log L/L)=O_\phi(L^{-1/2})$. This proves the bound for the terms with $e\ge2$. The estimates hold for complex vectors by replacing each coordinate by its absolute value.

Using
$$
\int_{\mathbb R}g(t)e^{-it\omega}dt
=\sqrt{2\pi}H\exp\left(-\frac{H^2\omega^2}{2}-\frac{3iT\omega}{2}\right),
$$
the lower bound for $|\log(q/q')|$ when $q\ne q'$, and Cauchy's inequality, we obtain, uniformly in $c$,
$$
\int_{\mathbb R}g(t)|R_T(t)|^2dt
=\sqrt{2\pi}H\|c\|^2(1+o(1)),
$$
with relative error at most $T e^{-L^{38}/8}$. By $|R_T|^2\le B\|c\|^2$ and Lemma~\ref{lem:logarithmic-boundary}, the normalised moments over $[T/2,5T/2]\setminus[T,2T]$ are at most $C_{\phi,j}T^{C_{\phi,j}}e^{-L^2/16}=o(1)$.

Choose $a=T/2+O(1)$ and $b=5T/2+O(1)$ outside radius-$T^{-2}$ neighbourhoods of shifted zeros. These excluded intervals have total length $O_\chi(T^{-2}\log T)<1$ in each unit interval. The logarithms on the vertical sides are then $O_{\chi,Y_0,\phi}(T\log^2T)$. We continue $|R_T(t)|^2$ analytically as the product of its two sums, keeping the conjugated coefficients constant. For $-1\le\operatorname{Im}t\le0$, this product satisfies
$\left|R_T(t)\sum_{q'}\overline{\chi(q')c_{q'}}q'^{it}\right| \le B^2\|c\|^2.$
On these sides, $g$ is $O(e^{-cL^2})$, so their contribution is $o(H\|c\|^2)$. By Cauchy's theorem and the boundary convergence in Lemma~\ref{lem:logarithmic-boundary}, we may move the integral to $\operatorname{Im}t=-1$.

On this line the Euler series converges absolutely. We expand $V_\phi^j$ and estimate its terms separately. Write $F_\chi(s)=\xi_\chi(s)/L(s,\chi)$ and
$$
G(t)=\pi^{-1}\left(\log F_\chi(1/2+i(t+h/2))-
\log F_\chi(1/2+i(t-h/2))\right).
$$
On $a\le\operatorname{Re}t\le b$, $|\operatorname{Im}t|\le H$, the completion factor $F_\chi$ has no zeros or poles in the corresponding $s$-region. By Stirling's formula in this sector, we have
$$
\begin{aligned}
\frac{F_\chi'}{F_\chi}(s)&=\frac12\log\frac Q\pi+
\frac12\frac{\Gamma'((s+a_\chi)/2)}{\Gamma((s+a_\chi)/2)}
+\boldsymbol1_{Q=1}\left(\frac1s+\frac1{s-1}\right)\\
&=\frac12\log\frac{Qs}{2\pi}+O_\chi(T^{-1}).
\end{aligned}
$$
Here $|s|\asymp T$ and $|\arg(s/2)|\le3\pi/4$. Thus $G$ is analytic throughout this rectangle and satisfies $G=i\phi+O_{\chi,\phi}(L^{-1})$ there. On $\operatorname{Im}t=-1$, we have
$$
V_\phi=G+\pi^{-1}\sum_{k\ge2}\chi(k)\beta(k)k^{-it},\qquad
\beta(k)=-\frac{2i\Lambda(k)}{\sqrt k\log k}\sin\frac{\pi\phi\log k}L.
$$
Let $\beta_r(k)$ be the coefficients of the $r$th power of $\sum_{k\ge2}\beta(k)k^{-it}$, with $\beta_0(1)=1$. Then
$$
\sum_k\frac{|\beta_r(k)|}k\le
\left(\frac{2\pi\phi}L\sum_{k\ge2}\frac{\Lambda(k)}{k^{3/2}}\right)^r
\le(C_\phi/L)^r.
$$
Consider a term with $r$ Euler-series factors and $m=j-r$ factors from $G$. If $\omega=\log(kq/q')\ne0$, then $|\omega|\ge1/(2B)$. We shift the integral of $gG^me^{-it\omega}$ to $\operatorname{Im}t=-\operatorname{sgn}(\omega)H$, where the exponential decreases. At this stage only $G$ and a single exponential are moved, after the Euler series has been expanded. The integral of $|g|$ along the horizontal segment is at most $e^{1/2}\sqrt{2\pi}H$. On the vertical sides, $|g|\le e^{-L^2/4}$ and $|e^{-it\omega}|\le e^{-\omega}$. Hence
$$
\begin{aligned}
\frac1{\sqrt{2\pi}H}\left|\int_{a-i}^{b-i}gG^me^{-it\omega}dt\right|
&\ll_{\phi,j}e^{-H|\omega|}+e^{-L^2/4}e^{-\omega},\\
e^{-H|\omega|}&\le e^{-L^{19}/4}\frac{q'}{kq}.
\end{aligned}
$$
Since $\sum_{q,q'}|c_qc_{q'}|q'/q\le B(\sum|c_q|)^2\le B^2\|c\|^2$, the oscillating terms and the integrals over their vertical sides contribute in total
$O_{\phi,j}\left(H\|c\|^2B^2(e^{-L^{19}/4}+e^{-L^2/4})\right) =o(H\|c\|^2).$
The remaining sum over $k$ is $\sum_k|\beta_r(k)|/k$, bounded above. This justifies shifting finite Euler sums and then passing to the full series.

For $\omega=0$, we have $q'=kq$. To identify these terms with $K^r$, we check the intermediate indices. If $q=m/n$ is reduced and $d\mid k$, the denominator of $dq$ is $n/(n,d)\le n$, and its numerator $md/(n,d)$ divides that of $kq$, prime by prime. Both therefore satisfy the size bounds and are coprime to $Q$. For integers, $dq\le kq\le N$. All intermediate indices are included, even where the coefficient vanishes.

By cancellation of the character factors, the sum of these ordered factorisations is $\sum_{q,q'}(K^r)_{q',q}c_q\overline{c_{q'}}$. Its sum of absolute values is at most
$$
\sum_{q,q'}(|K|^r)_{q',q}|c_qc_{q'}|
\le\||K|\|^r\|c\|^2=O_{\phi,j}(\|c\|^2).
$$ Since $(\sqrt{2\pi}H)^{-1}\int_{a-i}^{b-i}gG^m dt=(i\phi)^m+O_{\phi,j}(L^{-1}+e^{-L^2/4})$, replacing $G^m$ by $(i\phi)^m$ changes the integral by $O_{\phi,j}(H\|c\|^2/L)$. Binomial expansion and division by the mean-square integral of $R_T$ prove the assertion.

\end{proof}

\section{Evaluating the prime sums}\label{sec:integer-moments}

We choose coefficients supported on squarefree integers, so a nonzero term in a matrix power involves multiplication by distinct primes. Let $f$ be a nonzero continuous complex-valued function on $[0,1]$, and set
$$
c_n=\frac{\mu(n)}{\sqrt n}f\left(\frac{\log n}{\log N}\right)
\qquad(n\le N,\ (n,Q)=1).
$$
Lemma~\ref{lem:arithmetic-sums} evaluates $\sum_n|c_n|^2$. Multiplication by $p\nmid n$ changes $\log n/\log N$ by $u=\log p/\log N$ and reverses the M\"obius sign. The coefficient $-2i$ in $K$ therefore becomes $2i$. The limiting prime sum is the integral in
\begin{equation}\label{eq:A-phi}
(A_\phi f)(x)=i\phi f(x)+\frac{2i}\pi
\int_0^x\frac{\sin(\pi\phi u)}u f(x-u)\,du.
\end{equation}
The sine quotient is extended continuously at $u=0$, with value $\pi\phi$.

\begin{proposition}\label{prop:integer-limit}
For these coefficients and each fixed integer $j\ge0$,
$$
\sum_{\substack{n\le N\\(n,Q)=1}}|c_n|^2\sim\frac6{\pi^2}\prod_{p\mid Q}\frac p{p+1}\log N\int_0^1|f(x)|^2dx,
\qquad
\lim_{T\to\infty}\mathbb E_T[V_\phi^j]
=\frac{\displaystyle\int_0^1(A_\phi^jf)(x)\overline{f(x)}\,dx}{\|f\|_2^2}.
$$
The spectrum of $A_\phi$ is $\{i\phi\}$.
\end{proposition}
\begin{proof}
By Lemma~\ref{lem:arithmetic-sums} with $v=0$, we have the asserted asymptotic for $\sum_n|c_n|^2$. Proposition~\ref{prop:finite-moments} then reduces the moment formula to evaluating $\sum_{q,q'}(K^r)_{q',q}c_q\overline{c_{q'}}$ for each fixed $r$.

Since both $q$ and $q'$ are squarefree, the primes by which we multiply must be distinct and coprime to $q$. Thus, for $r\ge1$,
$$
\begin{aligned}
\sum_{q,q'}(K^r)_{q',q}c_q\overline{c_{q'}}
=(2i)^r\sum_{\substack{p_1,\ldots,p_r\ \mathrm{distinct}\\
p_i\nmid Q,\ P=p_1\cdots p_r\le N}}
&\prod_{i=1}^r\frac{\sin(\pi\phi\log p_i/\log T)}{p_i}\\
\times\sum_{\substack{n\le N/P\\(n,PQ)=1}}
&\frac{\mu(n)^2}{n}
 f\left(\frac{\log n}{\log N}\right)
\overline{f\left(\frac{\log(nP)}{\log N}\right)}.
\end{aligned}
$$
The primes are ordered, and $\mu(n)\mu(nP)=(-1)^r\mu(n)^2$ when $(n,P)=1$, so the factor is $(2i)^r$.

To evaluate the prime sums separately, we remove $(n,P)=1$ and distinctness with error $o(\log N)$. Set $s_p=\sin(\pi\phi\log p/\log T)$. By Mertens' estimate in Lemma~\ref{lem:arithmetic-sums} and $|s_p|\ll_\phi\log p/\log N$, we have $\sum_{p\le N}|s_p|/p=O_\phi(1)$. Also,
$$
\sum_{\substack{n\le Z\\p\mid n}}\frac{\mu(n)^2}{n}
\le\frac1p\sum_{m\le Z/p}\frac{\mu(m)^2}{m}
\ll\frac{1+\log Z}{p}
$$
when $Z\ge p$, and the sum is zero otherwise. Summing over $p_1,\ldots,p_r$ and bounding the remaining prime sums shows that omitting $(n,P)=1$ introduces an error of at most
$$
O_{Q,f,\phi,r}\left(\log N\sum_p\frac{|s_p|}{p^2}\right)=o(\log N).
$$
Allowing repeated primes then introduces an error of at most
$$
O_{Q,f,\phi,r}\left(\log N\sum_p\frac{|s_p|^2}{p^2}\right)=o(\log N)
\qquad(r\ge2).
$$
Both errors are $o(\log N)$ by dominated convergence, since $s_p\to0$ for each fixed $p$, $|s_p|\le1$, and $\sum_pp^{-2}<\infty$.

By partial summation of Mertens' estimate in Lemma~\ref{lem:arithmetic-sums}, we have, as $N\to\infty$,
$$
\frac1{\log N}\sum_{\substack{p\le N\\p\nmid Q}}
\frac{\log p}{p}F\left(\frac{\log p}{\log N}\right)
\longrightarrow\int_0^1F(u)\,du
$$
for every continuous $F$. Since $\log N/\log T\to1$, we may multiply $F(u)$ by
$$
\frac{\sin(\pi\phi u\log N/\log T)}u
$$
whose continuous extension at zero converges uniformly to $\sin(\pi\phi u)/u$. For $r$ primes, approximate continuous functions on $[0,1]^r$ by finite sums of products of one-variable functions, and apply the prime-sum limit in each variable. Using the first limit in Lemma~\ref{lem:arithmetic-sums} uniformly for $v=u_1+\cdots+u_r$, we obtain, as $T\to\infty$,
$$
\begin{aligned}
\frac{\sum_{q,q'}(K^r)_{q',q}c_q\overline{c_{q'}}}{\left(\frac6{\pi^2}\prod_{p\mid Q}\frac p{p+1}\right)\log N}
\longrightarrow (2i)^r
\int_{\substack{u_i\ge0\\\sum_i u_i\le1}}
&\prod_{i=1}^r\frac{\sin(\pi\phi u_i)}{u_i}\\
\times\int_0^{1-\sum_i u_i}
&f(x)\overline{f(x+\textstyle\sum_i u_i)}\,dx\,du_1\cdots du_r.
\end{aligned}
$$
The inner integral tends to zero as $\sum_i u_i\to1$, so its extension by zero is continuous. Substituting $y=x+\sum_i u_i$ and dividing by $\pi^r$ identifies the limit with $\int_0^1((A_\phi-i\phi)^rf)(x)\overline{f(x)}\,dx$, by repeated use of \eqref{eq:A-phi}. The asymptotic for $\sum_n|c_n|^2$ and binomial expansion in Proposition~\ref{prop:finite-moments} prove the moment formula.

Finally, $|2i\sin(\pi\phi u)/(\pi u)|\le2\phi$, so $r$ repeated integrations in \eqref{eq:A-phi} have convolution bounded by $(2\phi)^ru^{r-1}/(r-1)!$. Extending $f$ by zero outside $[0,1]$ and applying Young's convolution inequality, we obtain
$$
\|(A_\phi-i\phi)^r\|\le\frac{(2\phi)^r}{r!}.
$$
Thus, for every nonzero $z$, the series
$\sum_{r\ge0}(A_\phi-i\phi)^r/z^{r+1}$ converges in norm and inverts $z-(A_\phi-i\phi)$. On the other hand, $A_\phi-i\phi$ is not invertible, since otherwise
$1\le\|(A_\phi-i\phi)^{-1}\|^r(2\phi)^r/r!\to0$.
Thus the spectrum is $\{i\phi\}$. The same factorial estimate in the supremum norm proves uniform convergence of the exponential series on continuous functions.
\end{proof}

\section{Completing the small-gap proof}\label{sec:two-line-moments}

Under the gap assumption, the moments in Proposition~\ref{prop:integer-limit} belong to a positive measure on $\operatorname{Im}w=0,1$. We calculate its density and prove uniqueness. The contradiction follows by choosing $f$ so that the lower density is negative.

We write $\langle f,g\rangle=\int_0^1f(x)\overline{g(x)}dx$, over the shorter interval when appropriate. Thus $\langle Bf,f\rangle=\int(Bf)\overline f$. The adjoint satisfies $\langle Bf,g\rangle=\langle f,B^*g\rangle$.

The spectrum of a bounded operator $B$ consists of the $z$ for which $z-B$ has no bounded inverse. For $F$ analytic near this set, Cauchy's formula is
\begin{equation}\label{eq:operator-cauchy}
F(B)=\frac1{2\pi i}\int_\Gamma F(z)(z-B)^{-1}\,dz,
\end{equation}
where $\Gamma$ surrounds the spectrum counterclockwise. For $|z|>\|B\|$, expanding $(z-B)^{-1}=\sum_{j\ge0}B^j/z^{j+1}$ recovers powers and exponentials. The identity $(z-B)^{-1}(\zeta-B)^{-1}=((z-B)^{-1}-(\zeta-B)^{-1})/(\zeta-z)$ proves the rules for products and compositions. Contours may be deformed by Cauchy's theorem.

\begin{lemma}\label{lem:two-lines}
Let $B$ be a bounded operator on a complex Hilbert space, with spectrum in $0<\operatorname{Im}w<1$, and let $g$ be a unit vector. There is a real signed measure $\nu$ on $\operatorname{Im}w=0,1$ such that $\int w^j\,d\nu=\langle B^jg,g\rangle$ for every $j\ge0$. Its densities with respect to the real coordinate are
\begin{align}
p_0(x)&=e^{\pi x}\operatorname{Im}\left\langle
\left(e^{\pi x}-e^{\pi B}\right)^{-1}g,g\right\rangle,\label{eq:two-line-density}\\
p_1(x)&=e^{\pi x}\operatorname{Im}\left\langle
\left(-e^{\pi x}-e^{\pi B}\right)^{-1}g,g\right\rangle.\nonumber
\end{align}
Both densities are continuous. The total variation of $\nu$ has finite exponential moments in $|x|$ of every rate less than $\pi$. Any positive measure on these lines with finite absolute moments of all orders and the same moments equals $\nu$. 
\end{lemma}
\begin{proof}
The substitution $r=e^{\pi w}$ takes the lower line to the positive real axis and the upper line to the negative real axis.

Set $C=e^{\pi B}$. For $\lambda$ outside the exponential of the spectrum of $B$, apply \eqref{eq:operator-cauchy} to $(\lambda-e^{\pi z})^{-1}$. Multiplication by $\lambda-C$ identifies the integral as $(\lambda-C)^{-1}$. Thus the spectrum of $C$ lies in the upper half-plane. Since the principal logarithm satisfies $\operatorname{Log}(e^{\pi z})=\pi z$ near the spectrum of $B$, the same formula gives $\operatorname{Log}C=\pi B$. 
Define $d\Omega(r)=\pi^{-1}\operatorname{Im}\langle(r-C)^{-1}g,g\rangle dr$. Its density is bounded near zero because $C$ is invertible. At infinity, $\langle(r-C)^{-1}g,g\rangle=r^{-1}+O(r^{-2})$ has real leading term, so the density is $O(r^{-2})$. To calculate the moments, integrate
$$
(\operatorname{Log}\zeta/\pi)^j
\left((\zeta-C)^{-1}-(\zeta-C^*)^{-1}\right)
$$
around an upper half-annulus. The large arc contributes $O((1+\log R)^j/R)$. On the small arc the inverses are bounded, while its length times the logarithmic factor tends to zero. The spectra of $C$ and $C^*$ lie above and below the axis, respectively, so only the first inverse contributes. By \eqref{eq:operator-cauchy}, we have
$$
\int_{\mathbb R}(\operatorname{Log}(r+i0)/\pi)^j\,d\Omega(r)
=\langle B^jg,g\rangle.
$$
For $r>0$, let $x=\log r/\pi$, and for $r<0$, let $x+i=\log|r|/\pi+i$. Since $|dr|=\pi e^{\pi x}dx$, we obtain \eqref{eq:two-line-density}. The case $j=0$ gives mass one. The preceding bounds for $(r-C)^{-1}$ at zero and infinity give decay $O(e^{-\pi|x|})$.

For uniqueness, let $\nu_0',\nu_1'$ be the restrictions of a positive measure with these moments to the lower and upper lines, parametrised by the real coordinate. First compare coefficients in the following formal series. We justify the exponential integrals below.
$$
\begin{aligned}
\int e^{tx}d\nu_0'+e^{it}\int e^{tx}d\nu_1'&=\langle e^{tB}g,g\rangle,\\
\int e^{tx}d\nu_0'+e^{-it}\int e^{tx}d\nu_1'&=\langle e^{tB^*}g,g\rangle.
\end{aligned}
$$
The second identity follows by conjugating coefficients, since moments of the real coordinate are real and $\overline{\langle B^jg,g\rangle}=\langle(B^*)^jg,g\rangle$. Multiply respectively by $e^{-it/2}$ and $e^{it/2}$, and add. Formal division by $2\cos(t/2)$ is valid because its constant coefficient is $2$. Thus
$$
\sum_{j\ge0}\frac{t^j}{j!}\int x^j d(\nu_0'+\nu_1')
=\frac{\langle e^{t(B-i/2)}g,g\rangle+
\langle e^{t(B^*+i/2)}g,g\rangle}{2\cos(t/2)}.
$$
The right-hand side is analytic for $|t|<\pi$, so by Cauchy's estimate, we have $\int x^{2j}d(\nu_0'+\nu_1')\le C_R(2j)!R^{-2j}$ for $R<\pi$. Positivity and monotone convergence yield
$$
\int\cosh(ax)d(\nu_0'+\nu_1')\le C_R\sum_{j\ge0}(a/R)^{2j}<\infty
\quad(0<a<R<\pi).
$$
Since $e^{a|x|}\le2\cosh(ax)$, the integrals $\int e^{tx}d\nu_k'(x)$ are holomorphic on $|\operatorname{Re}t|<\pi$. By the identity theorem, the two exponential identities above hold throughout this strip. At $t=is$, writing $\widehat\nu_k'(s)=\int e^{isx}d\nu_k'(x)$, we have
$$
\begin{aligned}
\widehat\nu_0'(s)+e^{-s}\widehat\nu_1'(s)&=\langle e^{isB}g,g\rangle,\\
\widehat\nu_0'(s)+e^s\widehat\nu_1'(s)&=\langle e^{isB^*}g,g\rangle.
\end{aligned}
$$
These equations determine both Fourier transforms for $s\ne0$, and continuity determines their values at zero. The exponential decay of $p_0,p_1$ justifies the same calculation for the signed measure. By uniqueness of Fourier transforms, the measures are equal.
\end{proof}

\begin{lemma}\label{lem:two-line-limit}
Suppose that every interval of length $2\pi\phi/\log T$ in $[T-1,2T+1]$ contains at most one zero for all sufficiently large $T$. For each nonzero continuous $f$ on $[0,1]$, there is a positive measure on the two lines, with finite absolute moments of every order, whose moments are
$$
\frac{\displaystyle\int_0^1(A_\phi^jf)(x)\overline{f(x)}\,dx}
{\displaystyle\int_0^1|f(x)|^2\,dx},\qquad j\ge0.
$$
\end{lemma}
\begin{proof}
We must show that passage to a limiting measure preserves every moment. Write $V_\phi=X+iN_\phi$, with $X$ real and $N_\phi\in\{0,1\}$. Expanding the odd and even powers, we have
$$
\begin{aligned}
\operatorname{Im}\mathbb E_T[V_\phi^{2j+1}]
&=\sum_{r=0}^j(-1)^{j-r}\binom{2j+1}{2r}
\mathbb E_T[X^{2r}\boldsymbol1_{N_\phi=1}],\\
\operatorname{Re}\mathbb E_T[V_\phi^{2j}]
&=\mathbb E_T[X^{2j}]+
\sum_{r=0}^{j-1}(-1)^{j-r}\binom{2j}{2r}
\mathbb E_T[X^{2r}\boldsymbol1_{N_\phi=1}].
\end{aligned}
$$
The moments on the left converge by Proposition~\ref{prop:integer-limit}. Starting with $\mathbb E_TN_\phi$, the first identity proves convergence of $\mathbb E_T[X^{2j}\boldsymbol1_{N_\phi=1}]$ by induction, since its coefficient is $2j+1$. The second bounds $\mathbb E_TX^{2j}$. In particular, $\mathbb E_T|V_\phi|^2\le\mathbb E_TX^2+1=O(1)$, so the mass of $|V_\phi|>R$ is $O(R^{-2})$, uniformly in $T$. For fixed $j$, choose an integer $m$ with $2m>j$. Since $N_\phi=0$ or $1$,
$$
\begin{aligned}
\mathbb E_T|V_\phi|^{2m}&\le2^{m-1}(\mathbb E_T|X|^{2m}+1)=O_{\phi,f,m}(1),\\
\mathbb E_T\left[|V_\phi|^j\boldsymbol1_{|V_\phi|>R}\right]
&\le R^{j-2m}\mathbb E_T|V_\phi|^{2m}\longrightarrow0
\end{aligned}
$$
uniformly in large $T$ as $R\to\infty$. Successive subsequence selection makes the weighted proportions with $\operatorname{Re}V_\phi\le x$ converge on each line for every rational $x$. Their right-continuous limits define measures of total mass one, by the second-moment bound. Approximation by step functions whose intervals have endpoints of measure zero proves convergence of each power on bounded intervals. The tail estimate extends convergence to the whole line and bounds every absolute moment.
\end{proof}
\subsection{Proof of Theorem~\ref{thm:small}}\label{sec:small}

\begin{proof}
Suppose $\mu>1/2$. By the reflection argument in Section~\ref{sec:results}, all zeros above a fixed positive height $Y_0$ lie on the critical line. Choose $1/2<\phi<\phi'<\min(\mu,1)$. Eventually the gaps exceed $2\pi\phi'/\log\gamma_n$, and $\phi'/\log(2T+1)>\phi/\log T$. Hence every interval of length $2\pi\phi/\log T$ in $[T-1,2T+1]$ contains at most one zero. Lemma~\ref{lem:two-line-limit} therefore gives a positive measure with the required moments for each nonzero continuous $f$.

To choose $f$, we calculate the inverse in \eqref{eq:two-line-density} at $x=0$. Set $a=\pi\phi$ and $C=e^{\pi A_\phi}$. The Laplace transform of $2i\sin(au)/u$ is $2i\arctan(a/s)$, and $e^{2i\arctan(a/s)}=(s+ia)/(s-ia)=1+2ia/(s-ia)$.
For $s>0$, the integral $\int_0^\infty e^{-su}|\sin(au)|du/u$ is finite, so we may integrate the exponential convolution series term by term. By uniqueness of Laplace transforms, this series consists of the identity plus integration against $2ia e^{iau}\,du$. Since the repeated integrals use only $0\le u\le x$, on $[0,1]$ we have
\begin{equation}\label{eq:half-exponential}
(Cf)(x)=e^{ia}\left(f(x)+2ia\int_0^xe^{ia(x-y)}f(y)dy\right).
\end{equation}
To make $(1-C)^{-1}g$ explicit, let $h_0(x)=e^{-iax}$. By \eqref{eq:half-exponential}, $Ch_0=e^{ia}e^{iax}$. Thus $(1-C)^{-1}g=h_0$ when
$$
g=h_0-Ch_0=e^{-iax}-e^{ia}e^{iax}=-2ie^{ia/2}\sin(a(x+1/2)).
$$
The spectrum of $C$ is $\{e^{ia}\}$. Since $\pi/2<a<\pi$, $1-C$ has a bounded inverse and $g\ne0$. By Lemma~\ref{lem:two-lines}, applied to $g/\|g\|_2$, we have
\begin{align}
p_0(0)&=\frac{\displaystyle\operatorname{Im}\int_0^1((1-C)^{-1}g)(x)\overline{g(x)}\,dx}{\|g\|_2^2}\nonumber\\
&=\frac{\displaystyle\operatorname{Im}\int_0^1(Ch_0)(x)\overline{h_0(x)}\,dx}{\|g\|_2^2}
=\frac{\sin a\sin(2a)}{a\|g\|_2^2}<0.\label{eq:small-density}
\end{align}
Indeed the imaginary part of $\int_0^1h_0(x)\overline{h_0(x)-(Ch_0)(x)}\,dx$ equals that of $\int_0^1(Ch_0)(x)\overline{h_0(x)}\,dx=e^{ia}\int_0^1e^{2iax}dx=e^{2ia}\sin a/a$.

Lemma~\ref{lem:two-line-limit} with $f=g$ gives a positive measure with these moments. By uniqueness in Lemma~\ref{lem:two-lines}, its lower density is $p_0$. By \eqref{eq:small-density} and continuity, $p_0$ is negative on an interval. This contradiction proves $\mu\le1/2$.
\end{proof}

\section{A positive proportion of small spacings}\label{sec:positive-proportion}

To count small spacings, we need a weight with fourth moment $O(T)$. We use a shorter Dirichlet polynomial, omitting integers with small prime factors. We show that the sign in \eqref{eq:small-density} persists, then apply H\"older's inequality to count spacings.

Throughout this section we assume the Riemann Hypothesis for $\zeta(s)$, and set $L=\log T$, $h=2\pi\phi/L$ and $\chi=1$ in $V_\phi$. Products over $q\mid n$ run over prime divisors. Constants may depend on the fixed functions and parameters.

\begin{lemma}\label{lem:proportion-absolute-moments}
Assume the Riemann Hypothesis. For each fixed $p>1$ and $\phi>0$,
\begin{equation}\label{eq:proportion-V-moments}
\int_T^{2T}|V_\phi(t)|^pdt\ll_{\phi,p}T.
\end{equation}
Let $0<\eta<b$ and let $\rho$ be a smooth function with compact support in $(\eta,b)$. For each fixed $p\ge1$,
\begin{equation}\label{eq:proportion-prime-moments}
\int_T^{2T}\left|\sum_{q\ \mathrm{prime}}
\frac{\rho(\log q/L)}{q^{1/2+it}}\right|^pdt\ll_{\rho,p}T.
\end{equation}
If $0<\theta<1$ and $f$ is smooth with compact support on $\mathbb R$ and vanishes for $x\ge\theta$, then
\begin{equation}\label{eq:proportion-R-moments}
R_T(t)=\sum_n\frac{\mu(n)}{n^{1/2+it}}
f\left(\frac{\log n}{L}\right)
\prod_{q\mid n}\rho\left(\frac{\log q}{L}\right)
\quad\hbox{satisfies}\quad
\int_T^{2T}|R_T(t)|^pdt\ll_{f,\rho,p}T.
\end{equation}
\end{lemma}
\begin{proof}
Integrating Rodgers's tail estimate under the Riemann Hypothesis \cite[Theorem~1.1]{Rodgers} against $p x^{p-1}dx$ gives
\begin{equation}\label{eq:proportion-zero-moments}
\int_T^{2T}\left(\sum_\gamma
\frac1{1+(L(\gamma-t)/(2\pi))^2}\right)^pdt\ll_pT.
\end{equation}
The sum is over all non-trivial zero ordinates, positive and negative, counted with multiplicity. Every zero with $|\gamma-t|<A/L$ contributes at least $(1+(A/(2\pi))^2)^{-1}$, so the corresponding count has every fixed moment $O_{A,p}(T)$. A finite dyadic covering and comparison of logarithmic scales give the same bounds on $[cT,CT]$, for fixed $0<c<C$.

Equation~\eqref{eq:proportion-zero-moments} bounds the imaginary part of $V_\phi$. To bound its real part, pair the factors in the even Hadamard product for $\xi(1/2+it)$. The factors with zeros at $t=\gamma$, $0<\gamma<T/2$, and at $t=-\gamma$, $0<\gamma\le5T/2$, contribute $O(h/T)$ each. For $\gamma>5T/2$, the derivative of the paired logarithm is $O(T/\gamma^2)$. Partial summation of the zero count bounds their total by $O(h\log T)=O_\phi(1)$. Let $F_T(u)$ count zero ordinates in $[T/2,5T/2]\cap(u-h/2,u+h/2)$. Integrating the interval indicators gives
$$
\operatorname{Re}V_\phi(t)
=\frac1\pi\operatorname{p.v.}\int_{\mathbb R}\frac{F_T(u)}{t-u}\,du+O_\phi(1).
$$
The integral is a Cauchy principal value. Since $\int|F_T|^p\ll T$ by \eqref{eq:proportion-zero-moments}, the $L^p$ bound for the Hilbert transform \cite{Stein} proves \eqref{eq:proportion-V-moments} for $1<p<\infty$.

For \eqref{eq:proportion-prime-moments}, choose $\widehat F(u)=\rho(u)/u$, extended by zero to $u\le0$, with convention $\widehat F(u)=\int F(x)e^{-2\pi ixu}dx$. By linearity and the explicit formula \cite[pp.~108--109]{IK}\cite[Theorem~2.1]{Rodgers}, we have
$$
\sum_\gamma F\left(\frac{L(\gamma-t)}{2\pi}\right)
=C_F(t)-\frac1L\sum_{n\ge2}\frac{\Lambda(n)}{\sqrt n}
\widehat F\left(\frac{\log n}{L}\right)n^{-it}.
$$
Here $C_F$ is the gamma integral plus the pole terms
$F((\pm i/2-t)L/(2\pi))$. By Stirling's formula and the rapid decay of $F$, the gamma integral is $O_\rho(1)$. Repeated integration by parts bounds the pole terms by
$O_M(T^{b/2}(TL)^{-M}(1+L)^M)=O(1)$ for sufficiently large fixed $M$. The sum involving $\widehat F(-\log n/L)$ vanishes by its support. Prime terms have coefficient $\rho(\log q/L)$, and proper powers contribute at most
$$
\sum_{j\ge2}\frac1j\sum_{q\ \mathrm{prime}}
\frac{|\rho(j\log q/L)|}{q^{j/2}}\ll_\rho1.
$$
For $j=2$ use the prime harmonic sum on $T^{\eta/2}<q<T^{b/2}$, and for $j\ge3$ use absolute convergence. Since $|F(x)|\ll(1+x^2)^{-1}$, \eqref{eq:proportion-zero-moments} proves \eqref{eq:proportion-prime-moments}.

Replacing $\rho(u)$ by $\rho(u)e^{2\pi ivu}$ replaces $F(x)$ by $F(x+v)$. Integration by parts gives
$$
|F(x+v+iy)|\ll_M(1+|v|)^M(1+|x|)^{-M}
e^{2\pi b|y|}(1+|y|)^M.
$$
Applying the same estimates shows that the $L^p([T,2T])$ norm in \eqref{eq:proportion-prime-moments} is at most $T^{1/p}$ times a fixed power of $1+|v|$.

Each integer in \eqref{eq:proportion-R-moments} is below $T^\theta$ and has at most $k=\lfloor\theta/\eta\rfloor$ prime factors, all exceeding $T^\eta$. By Fourier inversion,
$$
R_T(t)=\sum_{r=0}^k(-1)^r\int_{\mathbb R}\widehat f(v)e_r(t,v)dv,
$$
where $e_r$ is the elementary symmetric polynomial of degree $r$ in
$\rho(\log q/L)q^{-1/2-it}e^{2\pi iv\log q/L}$. Higher degrees vanish after integration because $f$ vanishes beyond $\theta$. By Newton's identities, $e_r$ is a polynomial in the first $r$ power sums. The second is bounded absolutely by the prime harmonic estimate, as are higher powers by absolute convergence. Thus $|e_r|\ll_{\rho,r}(1+|S_1|)^r$, with $S_1$ the first power sum. Its modulated bound at order $pr$, Minkowski's inequality and the rapid decay of $\widehat f$ prove \eqref{eq:proportion-R-moments}.
\end{proof}

We next show that omitting primes below $T^\eta$ preserves the moment formula as $T\to\infty$ and then $\eta\downarrow0$. Let $0<\theta<1$, and define $A_\phi$ on $[0,\theta]$ by \eqref{eq:A-phi}. Choose smooth $0\le\psi,\omega\le1$, with $\psi=1$ on $[0,\theta]$, $\psi=0$ for $u\ge b>\theta$, $\omega=0$ for $u\le1$, and $\omega=1$ for $u\ge2$. Set
$$
\rho_\eta(u)=\psi(u)\omega(u/\eta),\qquad
c_n=\frac{\mu(n)}{\sqrt n}f\left(\frac{\log n}{L}\right)
\prod_{q\mid n}\rho_\eta\left(\frac{\log q}{L}\right).
$$
The support of $\rho_\eta$ lies in $(\eta/2,2b)$, as required by Lemma~\ref{lem:proportion-absolute-moments}. To average over length comparable to $T$, take
\begin{equation}\label{eq:proportion-averaged-Gaussian}
G_T(t)=\frac1{\sqrt{2\pi}H}\int_{5T/4}^{7T/4}
e^{-(t-u)^2/(2H^2)}du,\qquad H=T/L.
\end{equation}
Then $0\le G_T\le1$ and $\int_T^{2T}G_T(t)dt\sim T/2$.

\begin{lemma}\label{lem:proportion-rough-limit}
For each fixed polynomial $P$ and $f$ as in \eqref{eq:proportion-R-moments} with $\int_0^\theta|f(x)|^2dx>0$,
\begin{equation}\label{eq:proportion-rough-limit}
\lim_{\eta\downarrow0}\lim_{T\to\infty}
\frac{\int_T^{2T}G_T(t)|R_T(t)|^2P(V_\phi(t))dt}
{\int_T^{2T}G_T(t)|R_T(t)|^2dt}
=\frac{\int_0^\theta(P(A_\phi)f)(x)\overline{f(x)}dx}
{\int_0^\theta|f(x)|^2dx}.
\end{equation}
For every sufficiently small fixed $\eta$, the denominator divided by $T$ tends to a positive constant.
\end{lemma}
\begin{proof}
Average the weight in Proposition~\ref{prop:finite-moments} over centres $u\in[5T/4,7T/4]$, as in \eqref{eq:proportion-averaged-Gaussian}. The estimates there are uniform in $u$. The integrals outside $[T,2T]$ are $O(e^{-cL^2})$, since $u$ is at least $T/4$ from either endpoint. The lower bound for $|\log(kn/n')|$, when $kn\ne n'$, is unchanged, and the completion factor contributes $i\phi+O(L^{-1})$. Thus its moment formula holds with $G_T$, and
$$
\int_T^{2T}G_T(t)|R_T(t)|^2dt\sim\frac T2\sum_n|c_n|^2.
$$
The length $T^\theta<T/L^{20}$ is allowed for large $T$.

For fixed $\eta$, expand into the bounded number of prime factors and apply Mertens' estimate. The remaining cofactor is then averaged with respect to the finite measure $M_\eta$ defined by
$$
\int F(x)dM_\eta(x)=\sum_{k\ge0}\frac1{k!}
\int_{(0,b)^k}F(u_1+\cdots+u_k)
\prod_{i=1}^k\rho_\eta(u_i)^2\frac{du_i}{u_i}.
$$
The $k=0$ term is $F(0)$. Only finitely many terms contribute when $F$ is supported in $[0,\theta]$, and
$\sum_n|c_n|^2\longrightarrow\int|f|^2dM_\eta$. The terms obtained by multiplying successively by $r$ primes tend to
\begin{equation}\label{eq:proportion-inserted-primes}
(2i)^r\int_{\substack{u_1,\ldots,u_r\ge0\\v=u_1+\cdots+u_r\le\theta}}
\prod_{j=1}^r\frac{\rho_\eta(u_j)\sin(\pi\phi u_j)}{u_j}
\int_{[0,\theta-v]}f(x)\overline{f(x+v)}dM_\eta(x)
\,du_1\cdots du_r.
\end{equation}
The initial and final integers are squarefree, so the $r$ primes are distinct and occur to the first power. Requiring them to be coprime to the initial integer changes the sum by $O(\sum_{q>T^\eta}q^{-2})=o(1)$, since the other prime integrals are bounded. Multiplication by each prime reverses the M\"obius sign, replacing $-2i$ by $2i$. A prime introduced in this way contributes one factor $\rho_\eta$; a prime dividing the initial integer contributes two. The calculation in Proposition~\ref{prop:integer-limit} proves \eqref{eq:proportion-inserted-primes}.

To let $\eta$ decrease, normalise $M_\eta$ to mass one, obtaining $\nu_\eta$. Its Laplace transform is
$$
\int_0^\infty e^{-sx}d\nu_\eta(x)
=\exp\left(\int_0^b(e^{-su}-1)\rho_\eta(u)^2\frac{du}{u}\right)
\longrightarrow
\mathcal L(s)=\exp\left(\int_0^b(e^{-su}-1)\psi(u)^2\frac{du}{u}\right).
$$
The first moments are bounded by $\int_0^b\psi(u)^2du$, so subsequential limits have mass one. Each has Laplace transform $\mathcal L$, which determines it uniquely. Hence $\nu_\eta$ converges to a measure $\nu$ with this transform. Differentiating for $s>0$ and using uniqueness of the Laplace transform, we obtain, for $x>0$,
$$
x\,d\nu(x)=\left(\int_{[0,x]}\psi(x-y)^2d\nu(y)\right)dx.
$$
Moreover, $\nu(\{0\})=\lim_{s\to\infty}\mathcal L(s)=0$. The right side has a density, so $\nu$ has a density away from zero. Since $\psi=1$ on $[0,\theta]$, the function $F(x)=\nu([0,x])$ satisfies $xF'(x)=F(x)$ for $0<x<\theta$. Thus $F(x)=cx$ there, with $c>0$ because
$$
s\mathcal L(s)\longrightarrow
\frac{e^{-\gamma}}\theta
\exp\left(-\int_\theta^b\psi(u)^2\frac{du}{u}\right)>0,
$$
and the contribution of $x\ge\theta$ to $s\mathcal L(s)$ tends to zero. Here $\gamma$ is Euler's constant. The asymptotic follows by substituting $v=su$ in the integral over $(0,\theta)$, using
$\int_0^z(e^{-v}-1)dv/v=-\log z-\gamma+o(1)$.

Thus $\nu$ has no point masses and has constant positive density $c$ on $[0,\theta]$. In \eqref{eq:proportion-inserted-primes}, the cofactor integrals converge uniformly in $v$. Indeed, $f(x)\overline{f(x+v)}$, extended by zero beyond $\theta-v$, is a uniformly continuous family since $f(\theta)=0$. Approximate it by finitely many functions. The bound $|\sin(\pi\phi u)/u|\le\pi\phi$ permits dominated convergence in $u_1,\ldots,u_r$. The constant $c$ cancels on division by the denominator in \eqref{eq:proportion-rough-limit}. Substituting $y=x+u_1+\cdots+u_r$ gives $\int_0^\theta((A_\phi-i\phi)^rf)(y)\overline{f(y)}dy$. Binomial expansion proves \eqref{eq:proportion-rough-limit}, while $c\int_0^\theta|f|^2>0$ makes the denominator limit positive for small fixed $\eta$.
\end{proof}

\begin{proof}[Proof of Theorem~\ref{thm:positive-proportion}]
Choose $1/2<\phi<\min(1,1/2+\varepsilon)$ and $1/\phi-1<\theta<1$. We shall find a polynomial whose mean has negative real part, but whose real part is nonnegative whenever the interval contains at most one zero.

Formula \eqref{eq:half-exponential} holds on $[0,\theta]$, since its integrals use only points between zero and $x$. Set $a=\pi\phi$ and $g(x)=e^{-iax}-e^{ia}e^{iax}$. The same factorial bound gives spectrum $\{i\phi\}$. By Lemma~\ref{lem:two-lines}, the moments
$\int_0^\theta(A_\phi^jg)(x)\overline{g(x)}dx/\int_0^\theta|g(x)|^2dx$
have continuous real densities $p_0,p_1$ on the two lines, whose absolute values have finite exponential moments at every rate less than $\pi$. The calculation in \eqref{eq:small-density} gives
$$
p_0(0)=\frac{\sin(a\theta)\sin(a(1+\theta))}
{a\int_0^\theta|g(x)|^2dx}<0.
$$
The sign follows from $0<\phi\theta<1<\phi(1+\theta)<2$.

There is a real polynomial $q$ such that $\int q(x)^2p_0(x)dx<0$. To see this, polynomials are dense in $L^2(|p_0(x)|dx)$. A function orthogonal to every polynomial, multiplied by $|p_0|dx$, has an exponential moment by Cauchy's inequality. Its Fourier transform is analytic in a strip and all its derivatives vanish at zero, so uniqueness of Fourier transforms implies that the function is zero. Approximate a nonzero continuous function supported where $p_0<0$ by real polynomials. Their squares converge in $L^1(|p_0|dx)$, proving the claim.

Choose a real polynomial $b_0$ satisfying
$\operatorname{Im}b_0(x+i)=\operatorname{Re}(q(x+i)^2)$.
This is possible by solving from the highest degree down, since $\operatorname{Im}(x+i)^{m+1}$ has leading term $(m+1)x^m$. Then $P=q^2+ib_0$ satisfies
\begin{equation}\label{eq:proportion-polynomial-signs}
\operatorname{Re}P(x)=q(x)^2\ge0,\qquad
\operatorname{Re}P(x+i)=0,\qquad
\operatorname{Re}\int_0^\theta(P(A_\phi)g)(x)\overline{g(x)}dx<0.
\end{equation}
For the last inequality, integrate the first two identities against $p_0,p_1$ and use $\int q^2p_0<0$.

Approximate $g$ in $L^2([0,\theta])$ by a smooth, compactly supported $f$ vanishing for $x\ge\theta$. Boundedness of $P(A_\phi)$ preserves the last inequality in \eqref{eq:proportion-polynomial-signs} for sufficiently close $f$. By Lemma~\ref{lem:proportion-rough-limit}, choose $\eta>0$ small enough that the limiting quotient is negative. Its denominator is asymptotic to a positive multiple of $T$. Thus, for some fixed $\delta>0$ and all large $T$,
$$
\operatorname{Re}\int_T^{2T}G_T(t)|R_T(t)|^2P(V_\phi(t))dt
\le-\delta T.
$$
We fix $\phi,\theta,P,f,\eta$ in this order before letting $T\to\infty$.

Let $E_T$ be the set of $t\in[T,2T]$ whose interval $(t-h/2,t+h/2)$ contains at least two ordinates, counted with multiplicity. Outside $E_T$, the imaginary part of $V_\phi$ is zero or one, so by \eqref{eq:proportion-polynomial-signs}, we have $\operatorname{Re}P(V_\phi(t))\ge0$. By H\"older's inequality, Lemma~\ref{lem:proportion-absolute-moments}, and $G_T\le1$,
$$
\delta T\le\int_{E_T}G_T|R_T|^2|P(V_\phi)|dt
\le\left(\int_T^{2T}|R_T|^4dt\right)^{1/2}
\left(\int_T^{2T}|P(V_\phi)|^4dt\right)^{1/4}
\operatorname{meas}(E_T)^{1/4}
\ll T^{3/4}\operatorname{meas}(E_T)^{1/4}.
$$
Thus $\operatorname{meas}(E_T)\gg T$. Each such interval contains a consecutive pair, and the centres of intervals containing a fixed pair form a set of measure at most $h$. Hence there are $\gg T/h\gg T\log T$ spacings shorter than $h$. By the unit-interval zero count, requiring both ordinates in $[T,2T]$ removes only $O(\log T)$ pairs. For the remaining pairs,
$$
\frac{(\gamma_{n+1}-\gamma_n)\log\gamma_n}{2\pi}
<\phi\frac{\log(2T)}{\log T}<\frac12+\varepsilon
$$
for large $T$. This proves the theorem.
\end{proof}

\section{Large gaps}\label{sec:large}

Assume the Riemann Hypothesis for $L(s,\chi)$. If the weighted mean of the number of zeros in an interval is below one, some interval must be empty. We use symmetric coefficients indexed by reduced fractions $m/n$, with $m,n\le Y$. The completion factor contributes $\phi$ to this mean. The prime terms come from multiplying the numerator or cancelling a denominator factor, with equal contributions by symmetry.

Let $\ell>0$ weight each prime factor, and let $H$ be a real polynomial in five variables. Write $\Omega(n)$ for the number of prime factors of $n$, counted with multiplicity. For squarefree coprime $m,n\le Y$ with $(mn,Q)=1$, set
\begin{equation}\label{eq:large-marked-coefficients}
\begin{aligned}
c_{m/n}={}&\frac{\ell^{\Omega(mn)}}{\sqrt{mn}}\prod_{a=m,n}
H\left(\frac{\log a}{\log Y},
\sum_{p\mid a}\left(\frac{\log p}{\log a}\right)^{5/4},
\sum_{p\mid a}\left(\frac{\log p}{\log a}\right)^{3/2},\right.\\
&\hspace{1.3cm}\left.
\sum_{p\mid a}\left(\frac{\log p}{\log a}\right)^2,
\sum_{p\mid a}\left(\frac{\log p}{\log a}\right)^3\right).
\end{aligned}
\end{equation}
The product contains one value of $H$ for $m$ and one for $n$. Take all arguments zero at $a=1$, and other coefficients zero. The four prime sums lie in $[0,1]$ and describe how $\log a$ is divided among prime factors. They are near one when one prime accounts for most of $\log a$. With two nearly equal prime logarithms, they are near $2^{-1/4},2^{-1/2},1/2,1/4$. Their exponents were chosen by a finite search for a mean below one. Let $R_T(t)=\sum_q\chi(q)c_qq^{-it}$.

\subsection{Prime factors and the first moment}
For powers $j_1,\ldots,j_4$, let $(e_i)$ contain respectively $j_1,\ldots,j_4$ copies of $5/4,3/2,2,3$. On expanding the prime sums, indices assigned to the same prime form a block. This gives a partition $\mathcal P$ of the labelled indices.

\begin{lemma}\label{lem:large-marked-law}
There is a positive measure $P$ of total mass one on $[0,1]^4$ with moments
\begin{equation}\label{eq:large-partition-moments}
\begin{aligned}
&\int\prod_{i=1}^4t_i^{j_i}\,dP(t_1,t_2,t_3,t_4)\\
&\quad=\frac{\Gamma(\ell^2)}{\Gamma(\ell^2+5j_1/4+3j_2/2+2j_3+3j_4)}
\sum_{\mathcal P}\prod_{U\in\mathcal P}\ell^2\Gamma\left(\sum_{i\in U}e_i\right).
\end{aligned}
\end{equation}
All integrals with respect to $P$ are over $[0,1]^4$. The $j_i$ are nonnegative integers. When all $j_i=0$, the partition sum and product equal one. As $Y\to\infty$, the joint limiting measure of $v=\log n/\log Y$ and the four prime sums in \eqref{eq:large-marked-coefficients}, under the normalised weights $\ell^{2\Omega(n)}/n$ on squarefree $n\le Y$, is
$\ell^2v^{\ell^2-1}\,dv\,dP(t_1,t_2,t_3,t_4).$
Restricting to $(n,Q)=1$ leaves the limit unchanged. For the normalised pair weights $\ell^{2\Omega(mn)}/(mn)$ with $m,n\le Y$, $mn$ squarefree and $(mn,Q)=1$, the limiting averages factor into the product of the two one-variable averages.
\end{lemma}

\begin{proof}
By the Selberg--Delange method \cite[Chapter II.5]{Tenenbaum} and partial summation,
$$
\sum_{\substack{n\le x\\n\ \mathrm{squarefree}}}\frac{\ell^{2\Omega(n)}}n
=\frac{(\log x)^{\ell^2}}{\Gamma(\ell^2+1)}
\prod_p\left(1-\frac1p\right)^{\ell^2}\left(1+\frac{\ell^2}p\right)
+O_\ell((\log x)^{\ell^2-1}).
$$
The Euler product is positive. By monotonicity, the normalised sums up to $Y^v$ converge uniformly to $v^{\ell^2}$ on $[0,1]$. Restricting to integers coprime to fixed squarefree $d$ multiplies the constant by $\prod_{p\mid d}(1+\ell^2/p)^{-1}$.

Average a monomial with a continuous weight $\psi(\log n/\log Y)$. Replacing each denominator $\log n$ by $\log Y$ multiplies it by $v^{\sum_i e_i}$. A block $U$ assigned to a prime at $u=\log p/\log Y$ contributes $\ell^2u^{\sum_{i\in U}e_i-1}du$. Every block exponent exceeds one, so $\sum_{p\le x}(\log p)/p=\log x+O(1)$ bounds primes $p\le Y^\eta$ by $O_\ell(\eta)$. Above this cutoff, excluding repeated primes and requiring coprimality to the cofactor changes the sum by a relative $O_\ell(Y^{-\eta})$. Stieltjes integration and the uniform estimate for the sum over the remaining factor of $n$ give
$$
\frac{\Gamma(\ell^2)\prod_{U\in\mathcal P}\Gamma(\sum_{i\in U}e_i)}
{\Gamma(\ell^2+\sum_i e_i)}v^{\ell^2+\sum_i e_i-1}.
$$
After multiplication by $\ell^2$ from the cofactor density and from each block, and summation over partitions, the result is \eqref{eq:large-partition-moments} times
$\ell^2\int_0^1\psi(v)v^{\ell^2+\sum_i e_i-1}\,dv.$
On $v\ge\eta$, divide the integrands by $v^{\sum_i e_i}$. The mass with $v<\eta$ is $O_\ell(\eta^{\ell^2})+o(1)$, so we may let $\eta\downarrow0$. Polynomial approximation on the compact cube defines the measure and proves convergence along the full sequence. Its product with $\ell^2v^{\ell^2-1}dv$ gives the asserted joint limit.

For pairs, we insert $\sum_{d\mid(m,n)}\mu(d)$ and write $m=dm'$, $n=dn'$. Fixed $d$ changes sizes and sums by $o(1)$ away from zero, contributing
$$
\frac{\mu(d)\ell^{4\Omega(d)}}{d^2}\prod_{p\mid d}\left(1+\frac{\ell^2}p\right)^{-2}.
$$
The preceding estimate for $\sum_{n\le x}\ell^{2\Omega(n)}/n$ permits dominated convergence in $d$. After division by the sum of the weights, the positive, absolutely convergent product $\prod_p(1-\ell^4/(p+\ell^2)^2)$ cancels, leaving the product of the two one-variable averages. Excluding primes dividing $Q$ multiplies the single and pair constants by
$$
\prod_{p\mid Q}\left(1+\frac{\ell^2}p\right)^{-1}
\quad\hbox{and}\quad
\prod_{p\mid Q}\left(1+\frac{2\ell^2}p\right)^{-1}.
$$
These constants also cancel on division by the sum of the weights.
\end{proof}

To multiply the numerator by a prime, sum with weight $1/p$ over $p\nmid mn$, $mp\le Y$, and let $u=\log p/\log Y$. A continuous function bounded by $O(u)$ contributes $O(\eta)+o(1)$ from $p\le Y^\eta$. Above this cutoff, Mertens' theorem gives $du/u$, uniformly in the endpoint. Coprimality of $p$ to the cofactors changes the sum by a relative $O_\ell(1/p)$. Apply Lemma~\ref{lem:large-marked-law} and let $\eta\downarrow0$. Approximation by a finite set of continuous functions proves uniform convergence for uniformly bounded, equicontinuous families which are uniformly $O(u)$. The boundaries imposed by $m\le Y$ and $mp\le Y$ have measure zero, since the limiting measures for $\log m/\log Y$ and $\log(mp)/\log Y$ have densities. Each fixed $p\mid Q$ contributes $O_{Q,\ell}(1/\log Y)$, since the sine vanishes at this rate.

For logarithmic sizes $v=\log m/\log Y$ and $u=\log p/\log Y$, multiplication by $p$ changes the size to $v+u$. Dividing the prime logarithms by the new size changes the four sums to
$$
\begin{aligned}
t_1'&=\left(\frac v{v+u}\right)^{5/4}t_1+\left(\frac u{v+u}\right)^{5/4},&
t_2'&=\left(\frac v{v+u}\right)^{3/2}t_2+\left(\frac u{v+u}\right)^{3/2},\\
t_3'&=\left(\frac v{v+u}\right)^2t_3+\left(\frac u{v+u}\right)^2,&
t_4'&=\left(\frac v{v+u}\right)^3t_4+\left(\frac u{v+u}\right)^3.
\end{aligned}
$$
The mean square of the coefficients involves
$$
I=\int_0^1\int H(v,t_1,t_2,t_3,t_4)^2v^{\ell^2-1}\,dP(t_1,t_2,t_3,t_4)\,dv,
$$
and multiplication of the numerator by a prime involves
$$
\begin{aligned}
J(\phi)=\int_0^1\int_0^{1-v}\int{}&H(v,t_1,t_2,t_3,t_4)H(v+u,t_1',t_2',t_3',t_4')\\
&\times v^{\ell^2-1}\frac{\sin(\pi\phi u/2)}u\,dP(t_1,t_2,t_3,t_4)\,du\,dv.
\end{aligned}
$$
The product of the two $H$ factors and the sine is uniformly $O(u)$ and extends continuously by zero at $v=u=0$. The prime average just proved therefore applies.

\begin{proposition}\label{prop:large-marked-average}
If $I>0$, the coefficients \eqref{eq:large-marked-coefficients} satisfy
$$
\lim_{T\to\infty}\mathbb E_T[\operatorname{Im}V_\phi(t)]
=\phi-\frac{4\ell}{\pi}\frac{J(\phi)}I.
$$
\end{proposition}
\begin{proof}
By symmetry, it is sufficient to consider multiplication of the numerator by a prime. By Proposition~\ref{prop:finite-moments}, terms from $p^e$ with $e\ge2$ contribute $o(1)$ after division by $\sum_q|c_q|^2$, leaving
$$
-2i\sum_{\substack{m,n,p\colon\ p\ \mathrm{prime}\\mp,n\le Y\\mnp\ \mathrm{squarefree},\ (mnp,Q)=1}}
\frac{\sin(\pi\phi\log p/\log T)}{\sqrt p}\,c_{m/n}c_{mp/n}.
$$
Multiplication by $p$ multiplies $\ell^{\Omega(mn)}/\sqrt{mn}$ by $\ell/\sqrt p$. Together with $p^{-1/2}$ above, this gives $\ell/p$. Since $\log Y/\log T\to1/2$, the sine is $\sin(\pi\phi u/2)+o(u)$ uniformly. By Lemma~\ref{lem:large-marked-law}, the sum of squared coefficients is
$$
\sum_q|c_q|^2\sim
\frac{(\log Y)^{2\ell^2}I^2}{\Gamma(\ell^2)^2}
\prod_p\left(1-\frac1p\right)^{2\ell^2}\left(1+\frac{2\ell^2}p\right)
\prod_{p\mid Q}\left(1+\frac{2\ell^2}p\right)^{-1}.
$$
For the terms obtained by multiplying the numerator by a prime, $I^2$ is replaced by $-2i\ell J(\phi)I$, with the same arithmetic constant. For the terms obtained by cancelling a prime from the denominator, exchange the initial and final fractions, then interchange $m,n$, obtaining the same contribution. Their sum divided by $\sum_q|c_q|^2$ is $-4i\ell J(\phi)/I$. Its imaginary part divided by $\pi$ proves the formula.
\end{proof}

\subsection{Choice of the polynomial}
By Proposition~\ref{prop:large-marked-average}, it remains to choose $H$ with $I>0$ and
$(1-\phi)I+(4\ell/\pi)J(\phi)>0$. We use the polynomial space
$$
H(v,t_1,t_2,t_3,t_4)=\sum_{r,j_1,\ldots,j_4}b_{r,j_1,j_2,j_3,j_4}v^r\prod_{i=1}^4t_i^{j_i},
$$
where $0\le r\le12$, $j_i\ge0$ and $j_1+j_2+j_3+j_4\le5$. There are $\binom94=126$ monomials in the four prime sums and $13$ powers of $v$, hence $1638$ coefficients. Both $I$ and $J$ are quadratic forms in these coefficients. The summation indices below are integers, with $0\le r,s\le12$ and $0\le l_i\le k_i$. Both four-tuples $(j_1,\ldots,j_4),(k_1,\ldots,k_4)$ have nonnegative entries totalling at most $5$. By expansion, we have
$$
I=\sum_{r,s}\sum_{j_1,\ldots,j_4}\sum_{k_1,\ldots,k_4}\frac{b_{r,j_1,j_2,j_3,j_4}b_{s,k_1,k_2,k_3,k_4}}{\ell^2+r+s}
\int\prod_{i=1}^4t_i^{j_i+k_i}\,dP(t_1,t_2,t_3,t_4).
$$
By Taylor's theorem through degree $2K+2$, for integers $K\ge0$ we have
$$
\left|\sin x-\sum_{k=0}^K\frac{(-1)^kx^{2k+1}}{(2k+1)!}\right|
\le\frac{|x|^{2K+3}}{(2K+3)!}\qquad(x\in\mathbb R).
$$
Let $J_K$ be the integral defining $J$ with the sine replaced by this polynomial. In the second factor of $H$, choose $l_i$ factors from the original prime sums and $k_i-l_i$ from the prime by which we multiply. Set $v=xy$, $u=x(1-y)$, so that $x=v+u$ and $y=v/(v+u)$. The integrals in $x,y$ then separate, giving
$$
\begin{aligned}
J_K(\phi)={}&\sum_{r,s}\sum_{j_1,\ldots,j_4}\sum_{k_1,\ldots,k_4}\sum_{l_1,\ldots,l_4}
b_{r,j_1,j_2,j_3,j_4}b_{s,k_1,k_2,k_3,k_4}\\
&\times\prod_{i=1}^4\binom{k_i}{l_i}\int\prod_{i=1}^4t_i^{j_i+l_i}\,dP(t_1,t_2,t_3,t_4)\\
&\times\sum_{k=0}^K\frac{(-1)^k(\pi\phi/2)^{2k+1}}{(2k+1)!(\ell^2+r+s+2k+1)}\\
&\times\frac{\Gamma(\ell^2+r+5l_1/4+3l_2/2+2l_3+3l_4)}
{\Gamma(\ell^2+r+5k_1/4+3k_2/2+2k_3+3k_4+2k+1)}\\
&\times\Gamma\left(2k+\frac54(k_1-l_1)+\frac32(k_2-l_2)+2(k_3-l_3)+3(k_4-l_4)+1\right).
\end{aligned}
$$
The integral in $x$ is $(\ell^2+r+s+2k+1)^{-1}$. The gamma factors follow from the beta integral in $y$, and the moments are those in \eqref{eq:large-partition-moments}. Since arguments lie in $[0,1]$, $|H|\le\sum_{r,j_1,\ldots,j_4}|b_{r,j_1,j_2,j_3,j_4}|$ and
\begin{equation}\label{eq:marked-sine-remainder}
|J(\phi)-J_K(\phi)|\le
\frac{(\sum_{r,j_1,\ldots,j_4}|b_{r,j_1,j_2,j_3,j_4}|)^2}{\ell^2(\ell^2+1)}
\frac{(\pi\phi/2)^{2K+3}}{(2K+3)!}.
\end{equation}
Here $u^{2K+2}\le1$ and $\int_0^1v^{\ell^2-1}(1-v)dv=1/(\ell^2(\ell^2+1))$.

For fixed $\ell$, maximising $J_{32}/I$ is a finite generalised eigenvalue problem, with entries given by the formulas above. We take $\ell=39327/20000=1.96635$ and $\phi=3.63846723$. The value of $\ell$ was chosen numerically to increase this quotient. Its particular rational value is needed only for the finite calculation that follows.

Interval evaluation of these quadratic forms \cite{Arb}, at rational coefficients with denominators dividing $10^{80}$, gives $\sum|b_{r,j_1,j_2,j_3,j_4}|<7.46135$ and
$$
\begin{aligned}
\left|10^{34}I-1.000000000000000020968965813\right|&<10^{-27},\\
\left|10^{34}J_{32}-1.053854765588860867393555243\right|&<10^{-27}.
\end{aligned}
$$
In particular, $I>10^{-34}$ and
$$
1-\phi+\frac{4\ell}{\pi}\frac{J_{32}(\phi)}I
>2.15088826449\times10^{-9}.
$$
The subscript $32$ corresponds to $33$ odd terms, through degree $65$.
By \eqref{eq:marked-sine-remainder}, the remainder has power $67$ and satisfies
$|J-J_{32}|<4.273\times10^{-44}$.
Using the unrounded bound, its quotient contribution is less than
$1.069756\times10^{-9}$.
Subtracting this error leaves more than $1.08113226449\times10^{-9}$. Thus the explicit finite quadratic forms above satisfy
\begin{equation}\label{eq:marked-large-margin}
\sup_{\substack{H\text{ in the polynomial space above}\\I>0}}
\left(1-3.63846723+\frac{4\ell}\pi\frac{J(3.63846723)}I\right)
>1.0811\times10^{-9}.
\end{equation}

\begin{proof}[Proof of Theorem~\ref{thm:large}]
Fix $H$ satisfying \eqref{eq:marked-large-margin}. By Proposition~\ref{prop:large-marked-average}, the limiting weighted mean of the number of zeros is below one. By dominated convergence, $J$ is continuous in $\phi$, so the mean remains below one for some fixed $\phi'>3.63846723$. For large $T$, an interval of length $2\pi\phi'/\log T$ centred in $[T,2T]$ must therefore be empty. The zero-counting formula for fixed $\chi$ gives zeros in $[T/2,T-1]$ and $[2T+1,3T]$. The consecutive zeros on either side thus lie in $[T/2,3T]$, separated by at least $2\pi\phi'/\log T$. Since $\log\gamma_n/\log T\to1$ uniformly there, their normalised spacing is at least $\phi'+o(1)$. Hence $\lambda\ge\phi'>3.63846723$.
\end{proof}

\section*{Acknowledgements}
Both authors gratefully acknowledge support from the Heilbronn Institute for Mathematical Research. In this work, we used ChatGPT for the following tasks: as a sounding board during the development of the ideas, to proofread and probe the mathematical arguments, to run the optimisations found in the large gaps calculations, and to assist in the writing of and presentation of the paper. The authors take responsibility for the final contents of this work. 

Additionally, this paper will be split into two papers at a later date: the small gaps results will be merged with another paper on that problem, and the large gaps results will appear separately. Both will be rewritten by the authors. We leave this version publicly available while that work is carried out and will replace it when it is complete.


\begin{thebibliography}{99}
\raggedright

\bibitem{BaluyotAH}
S. A. C. Baluyot, \emph{On the pair correlation conjecture and the alternative hypothesis}, J. Number Theory \textbf{169} (2016), 183--226.

\bibitem{BennettZeros}
M. A. Bennett, G. Martin, K. O'Bryant and A. Rechnitzer, \emph{Counting zeros of Dirichlet $L$-functions}, Math. Comp. \textbf{90} (2021), 1455--1482.

\bibitem{BondarenkoHeapHalf}
A. Bondarenko and W. Heap, \emph{Half-sized gaps between zeros of the Riemann zeta function}, preprint, arXiv:2610.26038v1, 2026.

\bibitem{Bredberg}
J. Bredberg, \emph{Large gaps between consecutive zeros, on the critical line, of the Riemann zeta-function}, preprint, arXiv:1101.3197, 2011.

\bibitem{Bui11}
H. M. Bui, \emph{Large gaps between consecutive zeros of the Riemann zeta-function}, J. Number Theory \textbf{131} (2011), 67--95.

\bibitem{Bui14}
H. M. Bui, \emph{Large gaps between consecutive zeros of the Riemann zeta-function. II}, Acta Arith. \textbf{165} (2014), 101--122.

\bibitem{BGMM}
H. M. Bui, D. A. Goldston, M. B. Milinovich and H. L. Montgomery, \emph{Small gaps and small spacings between zeta zeros}, Acta Arith. \textbf{210} (2023), 133--153.

\bibitem{BM}
H. M. Bui and M. B. Milinovich, \emph{Gaps between zeros of the Riemann zeta-function}, Q. J. Math. \textbf{69} (2018), 403--423.

\bibitem{BMN}
H. M. Bui, M. B. Milinovich and N. Ng, \emph{A note on the gaps between consecutive zeros of the Riemann zeta-function}, Proc. Amer. Math. Soc. \textbf{138} (2010), 4167--4175.

\bibitem{CGGGH}
J. B. Conrey, A. Ghosh, D. Goldston, S. M. Gonek and D. R. Heath-Brown, \emph{On the distribution of gaps between zeros of the zeta-function}, Quart. J. Math. Oxford Ser. (2) \textbf{36} (1985), 43--51.

\bibitem{CGG84}
J. B. Conrey, A. Ghosh and S. M. Gonek, \emph{A note on gaps between zeros of the zeta function}, Bull. London Math. Soc. \textbf{16} (1984), 421--424.

\bibitem{CGG86}
J. B. Conrey, A. Ghosh and S. M. Gonek, \emph{Large gaps between zeros of the zeta-function}, Mathematika \textbf{33} (1986), 212--238.

\bibitem{CTB}
J. B. Conrey and C. L. Turnage-Butterbaugh, \emph{On $r$-gaps between zeros of the Riemann zeta-function}, Bull. London Math. Soc. \textbf{50} (2018), 349--356.

\bibitem{Davenport}
H. Davenport, \emph{Multiplicative Number Theory}, 3rd ed., revised by H.~L. Montgomery, Graduate Texts in Mathematics \textbf{74}, Springer-Verlag, New York, 2000.

\bibitem{FW12}
S. Feng and X. Wu, \emph{On gaps between zeros of the Riemann zeta-function}, J. Number Theory \textbf{132} (2012), 1385--1397.

\bibitem{FW13}
S. Feng and X. Wu, \emph{On large spacing of the zeros of the Riemann zeta-function}, J. Number Theory \textbf{133} (2013), 2538--2566.

\bibitem{GTT}
D. A. Goldston, T. S. Trudgian and C. L. Turnage-Butterbaugh, \emph{On the Montgomery--Odlyzko method regarding gaps between zeros of the zeta-function}, J. Math. Anal. Appl. \textbf{527} (2023), article 127548.

\bibitem{Hall}
R. R. Hall, \emph{A new unconditional result about large spaces between zeta zeros}, Mathematika \textbf{52} (2005), 101--113.

\bibitem{Inoue24}
S. Inoue, \emph{A note on $r$-gaps between zeros of the Riemann zeta-function}, Bull. London Math. Soc. \textbf{56} (2024), 2268--2277.

\bibitem{Inoue}
S. Inoue, \emph{Small gaps between consecutive zeros of the Riemann zeta-function}, preprint, arXiv:2604.05733, 2026.

\bibitem{IKT}
S. Inoue, H. Kobayashi and Y. Toma, \emph{Explicit extreme values of the argument of the Riemann zeta-function}, Math. Proc. Cambridge Philos. Soc. \textbf{181} (2026), 1149--1163.

\bibitem{IK}
H. Iwaniec and E. Kowalski, \emph{Analytic Number Theory}, American Mathematical Society Colloquium Publications \textbf{53}, American Mathematical Society, Providence, RI, 2004.

\bibitem{Arb}
F. Johansson, \emph{Arb: efficient arbitrary-precision midpoint-radius interval arithmetic}, IEEE Trans. Comput. \textbf{66} (2017), 1281--1292.

\bibitem{LagariasRodgers}
J. C. Lagarias and B. Rodgers, \emph{Higher correlations and the alternative hypothesis}, Q. J. Math. \textbf{71} (2020), 257--280.

\bibitem{Montgomery}
H. L. Montgomery, \emph{The pair correlation of zeros of the zeta function}, in \emph{Analytic Number Theory}, Proc. Sympos. Pure Math. \textbf{24}, American Mathematical Society, Providence, RI, 1973, 181--193.

\bibitem{MO}
H. L. Montgomery and A. M. Odlyzko, \emph{Gaps between zeros of the zeta function}, in \emph{Topics in Classical Number Theory}, Colloq. Math. Soc. J\'anos Bolyai \textbf{34}, North-Holland, 1984, 1079--1106.

\bibitem{Mueller}
J. Mueller, \emph{On the difference between consecutive zeros of the Riemann zeta function}, J. Number Theory \textbf{14} (1982), 327--331.

\bibitem{Ng}
N. Ng, \emph{Large gaps between the zeros of the Riemann zeta function}, J. Number Theory \textbf{128} (2008), 509--556.

\bibitem{Preobrazhenskii}
S. Preobrazhenski\u{\i}, \emph{A small improvement in the gaps between consecutive zeros of the Riemann zeta-function}, Res. Number Theory \textbf{2} (2016), article 28.

\bibitem{Rodgers}
B. Rodgers, \emph{Tail bounds for counts of zeros and eigenvalues, and an application to ratios}, Comment. Math. Helv. \textbf{92} (2017), 311--347.

\bibitem{Selberg}
A. Selberg, \emph{The zeta-function and the Riemann hypothesis}, in \emph{C. R. Dixi\`eme Congr\`es Math. Scandinaves 1946}, Jul. Gjellerups Forlag, Copenhagen, 1947, 187--200.

\bibitem{STTB}
A. Simoni\v{c}, T. S. Trudgian and C. L. Turnage-Butterbaugh, \emph{Some explicit and unconditional results on gaps between zeroes of the Riemann zeta-function}, Trans. Amer. Math. Soc. \textbf{375} (2022), 3239--3265.

\bibitem{Sou96}
K. Soundararajan, \emph{On the distribution of gaps between zeros of the Riemann zeta-function}, Quart. J. Math. Oxford Ser. (2) \textbf{47} (1996), 383--387.

\bibitem{Stein}
E. M. Stein, \emph{Singular Integrals and Differentiability Properties of Functions}, Princeton Mathematical Series \textbf{30}, Princeton University Press, Princeton, NJ, 1970.

\bibitem{Tenenbaum}
G. Tenenbaum, \emph{Introduction to Analytic and Probabilistic Number Theory}, 3rd ed., Graduate Studies in Mathematics \textbf{163}, American Mathematical Society, Providence, RI, 2015.

\bibitem{Titchmarsh}
E.~C. Titchmarsh, \emph{The Theory of the Riemann Zeta-function}, 2nd ed., edited by D.~R. Heath-Brown, Clarendon Press, Oxford, 1986.

\bibitem{Wu}
X. Wu, \emph{A note on the distribution of gaps between zeros of the Riemann zeta-function}, Proc. Amer. Math. Soc. \textbf{142} (2014), 851--857.

\end{thebibliography}
\end{document}